\documentclass[11pt]{article}
\usepackage{amsfonts}
\usepackage{amsmath}
\usepackage{amssymb}
\usepackage{xfrac}
\usepackage{setspace}
\numberwithin{equation}{section}
\newtheorem{Thm}{Theorem}[section]
\newtheorem{Cor}{Corollary}[section]
\newtheorem{Def}{Definition}[section]
\newtheorem{Lem}{Lemma}[section]
\newtheorem{Prop}{Proposition}[section]
\newtheorem{Remark}{Remark}[section]

\def\N{\mathbb{N}}

\def\R{\mathbb{R}}

\def\fd{\hfill$\square$}

\def\Rac{\mathcal{R}}
\newcommand{\prs}[2]{\mathop{\langle#1,#2\rangle}}
\date{}
\begin{document}
\title
{\emph{Green function and Poisson kernel associated to root systems  for annular regions}}
\author
{{\small\bf Chaabane REJEB}\footnote{Universit\'{e} de Tunis El Manar, Facult\'{e} des Sciences de Tunis, Laboratoire d'Analyse Math\'{e}matques et Applications LR11ES11, 2092 El Manar I, Tunis, TUNISIA, Laboratoire de Math\'ematiques et
Physique Th\'eorique CNRS-UMR 7350, Universit\'e de Tours,
Campus de Grandmont, 37200 Tours, FRANCE and Universit\'{e} de Sherbrooke, CANADA.
Email: chaabane.rejeb@gmail.com}}
\maketitle
\begin{abstract} Let $\Delta_k$ be the Dunkl Laplacian relative to a fixed root system $\Rac$ in $\R^d$, $d\geq2$, and to a nonnegative multiplicity function $k$ on $\Rac$. Our first purpose in this paper is to solve the $\Delta_k$-Dirichlet problem  for  annular regions. Secondly, we introduce and study the $\Delta_k$-Green function of the annulus and we prove that it  can be expressed  by means of $\Delta_k$-spherical harmonics. As applications, we obtain a Poisson-Jensen formula for $\Delta_k$-subharmonic functions and we study positive continuous solutions for a $\Delta_k$-semilinear problem.

\bigskip\noindent MSC (2020) primary: 31B05, 31B20, 31J05, 35J08; secondary: 31C45, 46F10, 47B39.

\end{abstract}
\section{Introduction}

Since the 90's, extensive  studies have been carried out on analysis  associated with Dunkl operators. These  are commuting differential-difference operators  on $\R^d$ introduced by C. F. Dunkl (see \cite{Dunkl1}).
 The Dunkl analysis includes especially   a generalization of the Fourier transform (called the Dunkl transform) and the Laplace operator known as the Dunkl Laplacian (and denoted by $\Delta_k$).\\
The Dunkl theory has many applications as well in mathematical physics and probability theory. In particular, it has been used  in the study of the Calogero-Moser-Sutherland and other integrable systems (see \cite{Diejen-Vinet, Etingof CMS}) and in the study of Markov processes generalizing Brownian motion (see \cite{Rosler-Voit}). \\
Recently,  a special interest has been  devoted to potential theory associated with the Dunkl Laplacian. The study focused on $\Delta_k$-harmonic functions (see \cite{Ben Chrouda, GalRej,GalRej2,Kods,Maslyouss, MajTrim}), on $\Delta_k$-Newton potential theory (including $\Delta_k$-subharmonic functions) (see \cite{GalRej3}) and on $\Delta_k$-Riesz potentials of Radon measures (see \cite{GalRejSifi}). More recently, by means of the $\Delta_k$-Newton kernel, the Green function of the open unit ball has been studied in \cite{Graczyk-Luks-Rosler}. Note that finding $\Delta_k$-Green functions for other  open sets is a rather difficult problem already in the case of the classical Laplace operator.  The aim of this paper is to show that we can determine the $\Delta_k$-Green function for annular regions in $\R^d$ by using $\Delta_k$-spherical harmonics as a crucial tool. \\
Let us assume throughout the paper that $d\geq2$. Let $A$ be the annulus
$$
A:=\{x\in\R^d,\ \rho<\|x\|<1\}\quad \text{with}\quad \rho\in (0,1).
$$
After giving some properties of  the $\Delta_k$-Green function $G_{k,A}$ of $A$, we will use it  to study   the  semilinear problem
$$
\left\{
  \begin{array}{ll}
    \Delta_k(u\omega_k)=\phi(.,u)\omega_k, & \hbox{in the sense of distributions} \\
    u=f, & \hbox{on}\ \partial{A},
  \end{array}
\right.
$$
where $\omega_k$ is a precise weight function (see  (\ref{weight function}) for its expression).\\
More precisely, under some assumptions on the function $\phi$, we will show that if $f\in \mathcal{C}(\partial{A})$ is nonnegative, this boundary problem has one and only one positive continuous solution on $A$  which  satisfies (see Theorem 5.2)
$$
 \forall\ x\in A,\quad u(x) + \int_AG_{k,A}(x,y)\phi(y,u(y))\omega_k(y)dy= P_{k,A}[f](x).
$$	
Here $P_{k,A}[f]$ is the unique solution in  $\mathcal{C}^{2}(A)\cap \mathcal{C}(\overline{A})$ of the boundary Dirichlet problem
$$
\left\{
  \begin{array}{ll}
    \Delta_ku=0, & \hbox{on}\ A,\\
    u=f, & \hbox{on}\ \partial{A},
  \end{array}
\right.
$$
that will be given explicitly in Section 3.
\medskip

This paper is organized as follows. In Section 2, we recall some  basics from Dunkl theory that will be used throughout the paper. In Section 3, we give an explicit solution of the boundary Dirichlet problem for the annulus. The Green function $G_{k,A}$ will be introduced and studied in Section 4. Some applications will be given in the last Section. Precisely, we will obtain a Poisson-Jensen  formula for $\Delta_k$-subharmonic functions in the annulus and we will study positive solutions of the above semilinear problem.

\section{Basics from Dunkl theory}
 We start by recalling   some useful facts in Dunkl theory. Let $\Rac$ be a root system in the Euclidian space $\R^d$, in the sense that $\Rac$ is a finite set in $\R^d\setminus\{0\}$ such that for every $\alpha\in \Rac$, $\Rac\cap\R\alpha=\{\pm\alpha\}$ and $\sigma_{\alpha}(\Rac)=\Rac$ (where $\sigma_{\alpha}$ is the reflection w.r.t. the hyperplane $H_{\alpha}$ orthogonal to $\alpha$). The subgroup $W\subset O(\R^d)$ generated by the reflections $\sigma_{\alpha}$, $\alpha\in\Rac$, is called the Coxeter-Weyl group associated to $\Rac$. We refer to (\cite{humph}) for more details on root systems and their Coxeter-Weyl groups.
\medskip

Let $k$ be a fixed nonnegative multiplicity function on $\Rac$ (i.e. $k$ is $W$-invariant). For $\xi\in\R^d$, the $\xi$-directional Dunkl operator associated to $(W,k)$ is defined by
$$
D_{\xi}f(x):=\partial_{\xi}f(x)+\sum_{\alpha\in\Rac_+}k(\alpha)\prs{\alpha}{\xi}\frac{f(x)-f(\sigma_\alpha.x)}{\prs{\alpha}{x}},\quad f\in \mathcal{C}^1(\R^d),
$$
where $\partial_{\xi}$ is the usual $\xi$-directional partial derivative and $\Rac_+$ is a positive subsystem.\\
Let us denote by $\mathcal{P}(\R^d)$ (resp. $\mathcal{P}_n(\R^d)$) the space of polynomial functions on $\R^d$ (resp. the space of homogeneous polynomials of degree $n\in\N$). \\
There exists a unique linear isomorphism $V_k$ from $\mathcal{P}(\R^d)$ onto itself such that
$V_k(\mathcal{P}_n(\R^d))=\mathcal{P}_n(\R^d)$ for every $n\in\N$, $V_k(1)=1$ and
\begin{equation}\label{Intertwining relation}
\forall\ \xi\in \R^d,\quad D_{\xi}V_k=\partial_{\xi}V_k.
\end{equation}
The operator $V_k$ is known as the Dunkl intertwining operator (see \cite{Dunkl2, DunklXu}). It has been extended to a topological isomorphism from $\mathcal{C}^{\infty}(\R^d)$ onto itself satisfying (\ref{Intertwining relation}) (see \cite{Trimeche1}). Furthermore, according to \cite{Rosler1}, for each $x\in\R^d$, there is a compactly supported probability measure $\mu_x$ on $\R^d$ such that
\begin{equation}\label{Integral rep Vk }
\forall\ f\in \mathcal{C}^{\infty}(\R^d),\quad V_k(f)(x)=\int_{\R^d}f(y)d\mu_x(y).
\end{equation}
If $W.x$ denotes the orbit of $x$ under the $W$-action and  $Co(x)$ its convex hull, then
\begin{equation}\label{support mu x}
\text{supp}\ \mu_x\subset Co(x)\subset \overline{B}(0,\|x\|).
\end{equation}

The Dunkl-Laplacian  is defined as $\Delta_k=\sum_{j=1}^dD_{e_j}^2$,
where $(e_j)_{1\leq j\leq d}$ is the canonical basis of $\R^d$. It can be expressed as follows
\begin{equation}\label{laplacian}
\Delta_kf(x)= \Delta f(x)+\sum_{\alpha\in R_+}k(\alpha)\Big(2\frac{\prs{\nabla f(x)}{\alpha}}{\prs{\alpha}{x}}-\|\alpha\|^2\frac{f(x)-f(\sigma_{\alpha}(x))}{\prs{\alpha}{x}^2}\Big),\quad f\in \mathcal{C}^2(\R^d),
\end{equation}
where $\Delta$ (resp. $\nabla$ ) is the usual Laplace (resp. gradient) operator (see \cite{Dunkl1, DunklXu}).
Note that if  $k$ is the zero function, the  Dunkl Laplacian reduces to the classical one which commutes with the action of $O(\R^d)$. For general $k\geq0$, $\Delta_k$ commutes with the $W$-action (see \cite{Rosler4}) i.e.
\begin{equation}\label{W-action}
\forall\ g\in W,\quad g\circ\Delta_k=\Delta_k\circ g.
\end{equation}

Let $L^2_k(S^{d-1})$, $d\geq 2$, be the Hilbert space endowed with the inner product
$$
 {\prs{p}{q}}_{k}:=\frac{1}{d_k}\int_{S^{d-1}}p(\xi)q(\xi)\omega_k(\xi)d\sigma(\xi).
 $$
We denote by $\|.\|_{L^2_k(S^{d-1})}$ the associated Euclidean norm. Here, $d\sigma$ is the surface measure on the unit sphere $S^{d-1}$, $\omega_k$ is the weight function given by
\begin{equation}\label{weight function}
\omega_k(x)=\textstyle\prod_{\alpha\in \Rac_+}|\prs{\alpha}{x}|^{2k(\alpha)}
\end{equation}
 and $d_k$ is the constant
\begin{equation}\label{dk}
d_k=\textstyle\int_{S^{d-1}}\omega_k(\xi)d\sigma(\xi).
\end{equation}
The function $\omega_k$ is $W$-invariant and homogeneous of degree $2\gamma:=2\textstyle\sum_{\alpha\in\Rac_+}k(\alpha)$.\\
Let us introduce the constant
\begin{equation}\label{lambda cte}
\lambda_k:=\frac{d}{2}+\gamma-1\geq0.
\end{equation}
Let $\mathcal{H}_{\Delta_k,n}(\R^d):=\mathcal{P}_n(\R^d)\cap Ker\Delta_k$ be the space of $\Delta_k$-harmonic polynomials, homogeneous  of degree $n$ on $\R^d$.
From \cite{DunklXu}, we know that if $n\neq m$, then $\mathcal{H}_{\Delta_k,n}(\R^d)\perp\mathcal{H}_{\Delta_k,m}(\R^d)$ in $L^2_k(S^{d-1})$. Moreover, for every $n\in\N$, we have
\begin{equation}\label{orthogonal decomp}
\mathcal{P}_n(\R^d)=\textstyle\bigoplus_{j=0}^{\lfloor n/2\rfloor}\|x\|^{2j}\mathcal{H}_{\Delta_k,n-2j}(\R^d).
\end{equation}
The restriction to the sphere $S^{d-1}$ of an element of $\mathcal{H}_{\Delta_k,n}(\R^d)$ is called a $\Delta_k$-spherical harmonic of degree $n$. The space of $\Delta_k$-spherical  harmonics of degree $n$ will be denoted by $\mathcal{H}_{\Delta_k,n}(S^{d-1})$. This space has a reproducing kernel $Z_{k,n}$ uniquely  determined  by the properties (see \cite{DaiXu, DunklXu})
  \begin{description}
    \item[i)] for each  $x\in S^{d-1}$, $Z_{k,n}(x,.)\in \mathcal{H}_{\Delta_k,n}(S^{d-1})$,
    \item[ii)] for every $f\in \mathcal{H}_{\Delta_k,n}(S^{d-1})$, we have
    \begin{equation}\label{reproducing kernel formula}
f(x)={\prs{f}{Z_{k,n}(x,.)}}_k=\frac{1}{d_k}\int_{S^{d-1}}f(\xi)Z_{k,n}(x,\xi)\omega_k(\xi)d\sigma(\xi),\quad x\in S^{d-1}.
\end{equation}
  \end{description}
  From this formula, we can see that
  \begin{equation}\label{Zn W-action}
  \forall\ g\in W,\quad \forall\ x,y\in S^{d-1},\quad Z_{k,n}(gx,gy)=Z_{k,n}(x,y).
  \end{equation}
In the classical case (i.e. $k=0$), $Z_{0,n}(x,.)$ is known as the zonal harmonic of degree $n$ (see \cite{Axler, DaiXu}).
Note that if $\{Y_{j,n}, j=1,\dots,h(n,d):=dim\mathcal{H}_{\Delta_k,n}(\R^d)\}$ is a real-orthonormal basis of $\mathcal{H}_{\Delta_k,n}(S^{d-1})$ in $L^2_k(S^{d-1})$, then
\begin{equation}\label{bon}
Z_{k,n}(x,y)=\sum_{j=1}^{h(n,d)}Y_{j,n}(x)Y_{j,n}(y).
\end{equation}
By means of the Dunkl intertwining operator and Gegenbauer polynomials, $Z_{k,n}$ is given explicitly  by (see \cite{DaiXu}, Theorem 7.2.6. or \cite{Xu97})
\begin{equation}\label{sphercial harrmonics Gegen}
\forall\ x,y\in S^{d-1},\quad Z_{k,n}(x,y)=\frac{(n+\lambda_k)(2\lambda_k)_n}{\lambda_k.n!}V_k\Big(P_{n}^{\lambda_k}\big(\prs{.}{y}\big)\Big)(x),
\end{equation}
where $\lambda_k$ is the constant given by (\ref{lambda cte}), $P_{n}^{\mu}$, $\mu>-1/2$, is the normalized  Gegenbauer polynomial (see \cite{DunklXu} p. 17) defined by
$$
P_{n}^{\mu}(x):=\frac{(-1)^n}{2^n(\mu+1/2)_n}(1-x^2)^{1/2-\mu}\frac{d^n}{dx^n}(1-x^2)^{n+\mu-1/2},
$$
and $(x)_n:=x(x+1)\dots(x+n-1)$ is the Pochhammer symbol.\\

At the end of this section, in order to simplify notations in the classical case $k=0$ we will write $L^2(S^{d-1})$ for $L^2_0(S^{d-1})$, $\mathcal{H}_{\Delta,n}$ for $\mathcal{H}_{\Delta_0,n}$, $|S^{d-1}|:=d_0$ the surface area of $S^{d-1}$ and $Z_n:=Z_{0,n}$.\\

\section{$\Delta_k$-Dirichlet problem for the annulus}
In this section, by introducing a Poisson type kernel,  we will solve the Dirichlet problem for the Dunkl Laplacian in annular regions
$$
A_{R_1,R_2}:=\{x\in\R^d:\ R_1<\|x\|<R_2\}.
$$
 Note that from the homogeneity property of $\Delta_k$:
$$
\delta_r\circ\Delta_k=r^{-2}\Delta_k\circ\delta_r,\quad \text{with}\quad \delta_r(f)(x):=f(rx),
$$
  it suffices to do this for  the annular region  $A=A_{\rho,1}$ with $\rho\in (0,1)$ fixed. \\
  Recall that  the $\Delta_k$-Poisson kernel of the unit ball (see \cite{DunklXu}) is given by
\begin{equation}\label{Poisson kernel unit ball} P_k(x,y)=\sum_{n=0}^{+\infty}Z_{k,n}(x,y)=\int_{\R^d}\frac{1-\|x\|^2}{\big(1-2\prs{x}{z}+\|x\|^2\big)^{\frac{d}{2}+\gamma}}d\mu_{y}(z),\quad (x,y)\in B\times {S^{d-1}}.
\end{equation}
From \cite{DunklXu}, we  know that
\begin{equation}\label{Poisson kernel as prob measure}
\frac{1}{d_k}\int_{S^{d-1}}P_k(x,\xi)\omega_k(\xi)d\sigma(\xi)=1.
\end{equation}
\medskip

We start by two preliminary useful  results.  For each $n\in \N$, the  restriction of the  Dunkl intertwining operator
$$
 V_k: \mathcal{H}_{\Delta,n}(\R^d)\longrightarrow \mathcal{H}_{\Delta_k,n}(\R^d)
$$
is a linear isomorphism.\\
In the first result, we will  estimate  the matrix-norms of this operator and of its inverse where the space  $\mathcal{H}_{\Delta,n}(\R^d)$ (resp. $\mathcal{H}_{\Delta_k,n}(\R^d)$) is endowed with the $L^2(S^{d-1})$-norm (resp. the $L_k^2(S^{d-1})$-norm).
 More precisely,
\begin{Prop} Let $n$ be a nonnegative integer.
\begin{description}
  \item[1.] For every $f\in\mathcal{H}_{\Delta,n}(\R^d)$, we have
  \begin{equation}\label{Vk estimate normeL2}
  \|V_k(f)\|_{L^2_k(S^{d-1})}\leq dim \mathcal{H}_{\Delta,n}(\R^d)\|f\|_{L^2(S^{d-1})}.
  \end{equation}
  \item[2.]  For every $f\in\mathcal{H}_{\Delta_k,n}(\R^d)$, we have
  \begin{equation}\label{Vk-1 estimate normeL2}
  \|V_k^{-1}(f)\|_{L^2(S^{d-1})}\leq \frac{(\gamma+\frac{d}{2})_n|S^{d-1}|}{(\frac{d}{2})_n} dim \mathcal{H}_{\Delta,n}(\R^d)\|f\|_{L^2_k(S^{d-1})}.
  \end{equation}
\end{description}
\end{Prop}
\emph{Proof:}\textbf{1)} Let $f\in \mathcal{H}_{\Delta,n}(\R^d)$. After rewriting the reproducing formula (\ref{reproducing kernel formula}) in the classical case (i.e. $k=0$),  applying it  to $f$ and using Fubini's theorem, we get
\begin{equation}\label{Vk-zonal}
V_k(f)(x)=\frac{1}{|S^{d-1}|}\int_{S^{d-1}} f(\xi)V_k[Z_n(.,\xi)](x)d\sigma(\xi),\quad x\in\R^d.
\end{equation}
But, from \cite{Axler}, Proposition 5.27, we have
$$
\forall\ z,\xi\in S^{d-1},\quad |Z_n(z,\xi)|\leq dim \mathcal{H}_{\Delta,n}(\R^d)
$$
 which implies that
\begin{equation}\label{estimate classical zonal}
\forall\ (z,\xi)\in\R^d\times S^{d-1},\quad |Z_{n}(z,\xi)|\leq \Big(dim \mathcal{H}_{\Delta,n}(\R^d)\Big)\|z\|^n.
\end{equation}
Thus, using the relations (\ref{Integral rep Vk }), (\ref{support mu x}), (\ref{Vk-zonal}) and (\ref{estimate classical zonal}) and the Cauchy-Schwarz inequality, we obtain
\begin{equation}\label{Vk ponctuelle estimate}
\forall\ x\in\R^d,\quad |V_k(f)(x)|\leq dim \mathcal{H}_{\Delta,n}(\R^d)\|f\|_{L^2(S^{d-1})}\|x\|^n.
\end{equation}
This implies that
$$
\|V_k(f)\|_{L^2_k(S^{d-1})}\leq dim \mathcal{H}_{\Delta,n}(\R^d)\|f\|_{L^2(S^{d-1})}.
$$
\textbf{2)} Let $f\in \mathcal{H}_{\Delta_k,n}(\R^d)$. By applying the classical case of the formula (\ref{reproducing kernel formula}) to $V_k^{-1}(f)$ and by using (\ref{estimate classical zonal}) and the Cauchy-Schwarz inequality, we deduce that
\begin{equation*}
\forall\ x\in \R^d,\quad |V_k^{-1}(f)(x)|\leq dim \mathcal{H}_{\Delta,n}(\R^d)\|V_k^{-1}(f)\|_{L^2(S^{d-1})}\|x\|^n.
\end{equation*}
Now,  using the following result (see \cite{DunklXu}, Proposition 5.2.8): for $p\in \mathcal{P}_n(\R^d)$ and $q\in \mathcal{H}_{\Delta,n}(\R^d)$,
then
$$
\frac{1}{|S^{d-1}|}\int_{S^{d-1}}p(\xi)q(\xi)d\sigma(\xi)=\frac{(\gamma+\frac{d}{2})_n|S^{d-1}|}{(\frac{d}{2})_nd_k}
\int_{S^{d-1}}p(\xi)V_k(q)(\xi)\omega_k(\xi)d\sigma(\xi)
$$
with $p=q=V_k^{-1}(f)$, we obtain
\begin{align*}
\|V_k^{-1}(f)\|_{L^2(S^{d-1})}^2&\leq \frac{(\gamma+\frac{d}{2})_n|S^{d-1}|}{(\frac{d}{2})_nd_k}\int_{S^{d-1}}|V_k^{-1}(f)(\xi)f(\xi)|\omega_k(\xi)d\sigma(\xi)\\
&\leq \frac{(\gamma+\frac{d}{2})_n|S^{d-1}|}{(\frac{d}{2})_nd_k} dim \mathcal{H}_{\Delta,n}(\R^d)\|V_k^{-1}(f)\|_{L^2(S^{d-1})}\int_{S^{d-1}}|f(\xi)|\omega_k(\xi)d\sigma(\xi)\\
&\leq \frac{(\gamma+\frac{d}{2})_n|S^{d-1}|}{(\frac{d}{2})_n} dim \mathcal{H}_{\Delta,n}(\R^d)\|V_k^{-1}(f)\|_{L^2(S^{d-1})}
\|f\|_{L^2_k(S^{d-1})}.
\end{align*}
This proves the desired relation.
\fd
\begin{Cor} The following inequality holds:
\begin{equation}\label{dimension inequality}
 \forall\ x,y\in S^{d-1},\quad |Z_{k,n}(x,y)|\leq \Big(\frac{(\gamma+\frac{d}{2})_n|S^{d-1}|}{(\frac{d}{2})_n}\Big)^2 \Big(dim \mathcal{H}_{\Delta,n}(\R^d)\Big)^5.
\end{equation}
\end{Cor}
\emph{Proof:} Let $\{Y_{j,n}\}_j$, $j=1,\dots,h(n,d)=dim\mathcal{H}_{\Delta_k,n}(\R^d)$, be a real-orthonormal basis of $\mathcal{H}_{\Delta_k,n}(\R^d)$ in $L^2_k(S^{d-1})$. Using (\ref{Vk ponctuelle estimate}) with $f=V_k^{-1}(Y_{j,n})$ and (\ref{Vk-1 estimate normeL2}), we deduce that
\begin{align*}
\forall\ x\in S^{d-1}, |Y_{j,n}(x)|&\leq dim \mathcal{H}_{\Delta,n}(\R^d)\|V_k^{-1}(Y_{j,n})\|_{L^2(S^{d-1})}\\
&\leq \frac{(\gamma+\frac{d}{2})_n|S^{d-1}|}{(\frac{d}{2})_n} \Big(dim \mathcal{H}_{\Delta,n}(\R^d)\Big)^2.
\end{align*}
Consequently, we obtain the result from (\ref{bon}).
\fd

\bigskip

\noindent Following the classical case $k=0$ (see \cite{Axler}), we define  the kernel $\mathbf{P}_{k,1}(.,.)$ on $A\times{S^{d-1}}$ by
\begin{equation}\label{Poisson kernel part1}
\mathbf{P}_{k,1}(x,\xi):=\sum_{n=0}^{+\infty}a_{k,n}(x)Z_{k,n}(x,\xi), \quad \text{with}\quad a_{k,n}(x)=\frac{1-\big(\frac{\|x\|}{\rho}\big)^{-2\lambda_k-2n}}{1-\rho^{2\lambda_k+2n}}.
\end{equation}
\begin{Prop}\label{Poisson kernel part 1 lemma} The kernel $\mathbf{P}_{k,1}$ satisfies the following properties
\begin{description}
  \item[i)] For each $\xi\in S^{d-1}$, $\mathbf{P}_{k,1}(.,\xi)$ is a  $\Delta_k$-harmonic function on $A$ and $\mathbf{P}_{k,1}(.,\xi)=0$  on $S(0,\rho)$.
  \item[ii)] For every $x\in A$ and $\xi\in S^{d-1}$,
  \begin{equation}\label{norme1 Poisson kernel annulus}
0\leq \mathbf{P}_{k,1}(x,\xi)\leq P_k(x,\xi).
\end{equation}
  \item[iii)] Let $x\in A$ and $\xi\in S^{d-1}$ fixed. Then
  \begin{equation}\label{W-action kernel 1}
  \forall\ g\in W,\quad \mathbf{P}_{k,1}(gx,g\xi)=\mathbf{P}_{k,1}(x,\xi).
  \end{equation}
\end{description}
\end{Prop}
\noindent\emph{Proof:}
 \textbf{i)} Clearly $\mathbf{P}_{k,1}(.,\xi)=0$ on $S(0,\rho)$.
On the other hand, for any $(x,\xi)\in A\times S^{d-1}$ we can write
$$
a_{k,n}(x)Z_{k,n}(x,\xi)=c_{1,n}Z_{k,n}(x,\xi)-c_{2,n}K_k[Z_{k,n}(.,\xi)](x),
$$
where $c_{1,n},c_{2,n}$ are two nonnegative constants and $K_k$ is the $\Delta_k$-Kelvin transform (see \cite{El Kamel}) given by
\begin{equation}\label{Kelvin transf def}
K_k[f](x)=\|x\|^{-2\lambda_k}f(x/{\|x\|^2})=\|x\|^{2-2\gamma-d}f(x/{\|x\|^2})
\end{equation}
and $f$ is a function defined on $\R^d\setminus\{0\}$. As the $\Delta_k$-Kelvin transform preserves the $\Delta_k$-harmonic functions  on $\R^d\setminus\{0\}$  (see \cite{El Kamel}), we deduce that the function $x\mapsto a_{k,n}(x)Z_{k,n}(x,\xi)$ is $\Delta_k$-harmonic on $A$.\\
According to \cite{DunklXu} (see also \cite{Axler} and \cite{DaiXu}), we know that
$$
dim\mathcal{H}_{\Delta,n}(\R^d)=dim\mathcal{H}_{\Delta_k,n}(\R^d)={n+d-1\choose{n}}-{n+d-3\choose{n-2}}.
$$
Hence, we have
$$
\lim_{n\rightarrow+\infty}n^{2-d}dim\mathcal{H}_{\Delta,n}(\R^d)=\frac{2}{(d-2)!}.
$$
Moreover, we have
$$
\lim_{n\rightarrow+\infty}n^{-\gamma}\frac{(\gamma+\frac{d}{2})_n}{(\frac{d}{2})_n}=\lim_{n\rightarrow+\infty}n^{-\gamma}\frac{\Gamma(d/2)}{\Gamma(\gamma+d/2)}
\frac{\Gamma(d/2+\gamma+n)}{\Gamma(d/2+n)}=\frac{\Gamma(d/2)}{\Gamma(\gamma+d/2)}.
$$
Consequently, from (\ref{dimension inequality}), there exists $C=C(d,\gamma)>0$ such that
\begin{equation}\label{k-sphercial harmonic ineq}
\forall \ x\in\R^d,\quad\forall\ y\in S^{d-1}, \quad |Z_{k,n}(x,y)|\leq C n^{5d+2\gamma-10}\|x\|^n.
\end{equation}
This inequality as well as the fact that $0\leq a_{k,n}(x)<1$ imply that the series
$$
\textstyle\sum_{n\geq0}a_{k,n}(x)Z_{k,n}(x,\xi)
$$
 converges  uniformly on $\overline{A}_{\rho,R}\times S^{d-1}$ for every $R\in (\rho,1)$.
Then , by Corollary 3.3 in \cite{GalRej}, the function $\mathbf{P}_{k,1}(.,\xi)$ is $\Delta_k$-harmonic on $A$.\\
\textbf{ii)} For $\varepsilon>0$ small enough and $\xi\in S^{d-1}$, consider the function
$$
h_{\varepsilon}(x):=\textstyle\sum_{n\geq0}a_{k,n}(x)Z_{k,n}((1-\varepsilon)x,\xi).
$$
As above, from the inequality (\ref{k-sphercial harmonic ineq}) and  the homogeneity of $Z_{k,n}(.,\xi)$, we see that $h_{\varepsilon}$ defines a $\Delta_k$-harmonic function in the annular region $A_{\rho,R}$ with $R=(1-\varepsilon)^{-1}$. Furthermore,  $h_\varepsilon=0$ on $S(0,\rho)$ and if $x\in S^{d-1}$, then
\begin{equation}\label{relation P1-P}
h_{\varepsilon}(x)=\textstyle\sum_{n\geq0}Z_{k,n}((1-\varepsilon)x,\xi)=P_k((1-\varepsilon)x,\xi).
\end{equation}
where $P_k$ is the $\Delta_k$-Poisson kernel of the unit ball (see \cite{DunklXu}). In particular, $h_{\varepsilon}\geq 0$ on $S^{d-1}$. Consequently, by
 the weak minimum principle for $\Delta_k$-harmonic functions (see \cite{GalRej} or \cite{Rosler3}), we deduce that
$$
\forall\ x\in A,\quad h_{\varepsilon}(x)\geq 0.
$$
On the other hand, for each fixed $(x,\xi)$ in $A\times S^{d-1}$, we have
\begin{align*}
|\mathbf{P}_{k,1}(x,\xi)-h_{\varepsilon}(x)|&\leq \sum_{n\geq1}(1-(1-\varepsilon)^n)a_{k,n}(x)|Z_{k,n}(x,\xi)|\\
&\leq C\sum_{n\geq1}(1-(1-\varepsilon)^n)n^{5d+2\gamma-10}\|x\|^n.
\end{align*}
Hence, by the monotone convergence theorem we have $\mathbf{P}_{k,1}(x,\xi)=\lim_{\varepsilon\rightarrow0}h_{\varepsilon}(x)$. Finally,  we obtain $\mathbf{P}_{k,1}\geq0$ on $A\times S^{d-1}$.\\
$\bullet$ For $\xi\in S^{d-1}$ fixed, the function $x\mapsto P_k((1-\varepsilon)x,\xi)-h_{\varepsilon}(x)$ is $\Delta_k$-harmonic on $A$. Moreover, since $P_k$ is a nonnegative kernel, we have
$$
\forall\ x\in S(0,\rho),\quad P_k((1-\varepsilon)x,\xi)- h_{\varepsilon}(x)=P_k((1-\varepsilon)x,\xi)\geq0.
$$
By (\ref{relation P1-P}), $x\mapsto P_k((1-\varepsilon)x,\xi)- h_{\varepsilon}(x)$ is the zero function on $S^{d-1}$.
So, the weak maximum principle implies that
$$
\forall \ x\in A,\quad P_k((1-\varepsilon)x,\xi)\geq h_{\varepsilon}(x).
$$
Letting $\varepsilon\longrightarrow0$, we obtain $P_k(.,\xi)\geq \mathbf{P}_{k,1}(.,\xi)$ on $A$.
\medskip

\noindent
\textbf{iii)} The result follows immediately from (\ref{Zn W-action}).
\fd

\begin{Prop}\label{Poisson integral part 1} Let $f$ be a continuous function on $S^{d-1}$. Then the function
\begin{equation}\label{solution1}
\mathbf{P}_{k,1}[f](x)=\frac{1}{d_k}\int_{S^{d-1}}\mathbf{P}_{k,1}(x,\xi)f(\xi)\omega_k(\xi)d\sigma(\xi)
\end{equation}
is the unique solution in $\mathcal{C}^2(A)\cap\mathcal{C}(\overline{A})$ of the boundary Dirichlet problem
$$
\left\{
  \begin{array}{ll}
    \Delta_ku=0, & \hbox{on}\ A; \\
     u=f, & \hbox{on}\ S^{d-1}\\
     u=0, & \hbox{on}\ S(0,\rho).
  \end{array}
\right.
$$
\end{Prop}
 \emph{Proof:}
 The uniqueness follows from the weak maximum principle for $\Delta_k$-harmonic functions (see \cite{GalRej} or \cite{Rosler3}). The inequality (\ref{k-sphercial harmonic ineq})  allowed us to write for any $x\in A$ that
$$
\mathbf{P}_{k,1}[f](x)=\sum_{n=0}^{+\infty}u_n(x),\quad \text{with}\quad u_n(x)=\frac{a_{k,n}(x)}{d_k}\int_{S^{d-1}}Z_{k,n}(x,\xi)f(\xi)\omega_k(\xi)d\sigma(\xi).
$$
By differentiation theorem under integral sign, the functions $u_n$ are $\Delta_k$-harmonic on $A$. Moreover, by (\ref{k-sphercial harmonic ineq}) we have
$$
\forall\ n,\quad |u_n(x)|\leq C\|f\|_{\infty} n^{5d+2\gamma-10}\|x\|^n.
$$
This proves that the series $\sum_{n\geq0}u_n$ converges uniformly on each closed annular region $\overline{A}_{\rho,R}$ whenever $R\in (\rho,1)$. Then, we conclude that $\mathbf{P}_{k,1}[f]$ is $\Delta_k$-harmonic on $A$.\\
On the other hand, it is easy to see that $\mathbf{P}_{k,1}[f]=0$ on $S(0,\rho)$.\\
It remains to prove that for every $\xi\in S^{d-1}$, $\lim_{x\rightarrow \xi}\mathbf{P}_{k,1}[f](x)=f(\xi)$.\\
- If $f\in \mathcal{H}_{\Delta_k,m}(\R^d)$, then $u_n=0$ if $n\neq m$ and $u_m(x)=a_{k,m}(x)f(x)=\mathbf{P}_{k,1}[f](x)$. Therefore, $\mathbf{P}_{k,1}[f]=f$ on $S^{d-1}$.\\
- If $f\in \mathcal{P}_m(\R^d)$, then by (\ref{orthogonal decomp}), there exist $f_1,\dots,f_m$, with $f_j\in \mathcal{H}_{\Delta_k,n-2j}(\R^d)$ such that
$$
f(x)=\textstyle\sum_{j=0}^{[m/2]}\|x\|^{2j}f_j(x).
$$
This implies  that $\mathbf{P}_{k,1}[f]=f$ on $S^{d-1}$. \\
- If $f$ is an arbitrary  polynomial function, the result also holds. \\
- Suppose that $f$ is  a continuous function on  $S^{d-1}$ and let $p$ be a polynomial function. \\
By (\ref{norme1 Poisson kernel annulus}) and (\ref{Poisson kernel as prob measure}) we have
\begin{align*}
|\mathbf{P}_{k,1}[f](x)-f(x)|&\leq |\mathbf{P}_{k,1}[f](x)-\mathbf{P}_{k,1}[p](x)|+ |\mathbf{P}_{k,1}[p](x)-p(x)|+|p(x)-f(x)|\\
&\leq 2\|f-p\|_{\infty}+|\mathbf{P}_{k,1}[p](x)-p(x)|.
\end{align*}
This inequality as well as the Stone-Weierstrass theorem show that  $\lim_{x\rightarrow \xi}\mathbf{P}_{k,1}[f](x)=f(\xi)$ for every $\xi\in S^{d-1}$.
This completes the proof.
\fd

\bigskip

Now, for $x\in A$ and $\xi\in S^{d-1}$, consider the functions
$$
 b_{k,n}(x)=\|x\|^{-n}\big(\frac{\|x\|}{\rho}\big)^{-2\lambda_k-n}\frac{1-\|x\|^{2\lambda_k+2n}}{1-\rho^{2\lambda_k+2n}}=\rho^{-n}(1-a_{k,n}(x))
$$
 and
\begin{equation}\label{Poisson kernel part 2}
\mathbf{P}_{k,2}(x,\xi):=\sum_{n=0}^{+\infty}b_{k,n}(x)Z_{k,n}(x,\xi).
\end{equation}
By means of the Poisson kernel of the unit ball, we can  write
\begin{equation}\label{link with Poisson kernel ball}
\mathbf{P}_{k,2}(x,\rho\xi)=P_k(x,\xi)-\mathbf{P}_{k,1}(x,\xi).
\end{equation}
This relation as well as the properties of $P_k$ and $\mathbf{P}_{k,1}$ prove that $\mathbf{P}_{k,2}(.,\rho\xi)$ is a nonnegative $\Delta_k$-harmonic function on $A$ with $\mathbf{P}_{k,2}(.,\rho\xi)=0$ on $S^{d-1}$.

\medskip
Let $f$ be  a continuous function on $S(0,\rho)$ and define the function
\begin{equation}\label{solution2}
\mathbf{P}_{k,2}[f](x)=\frac{1}{d_k}\int_{S^{d-1}}\mathbf{P}_{k,2}(x,\rho\xi)f(\rho\xi)\omega_k(\xi)d\sigma(\xi), \quad x\in A.
\end{equation}
Using (\ref{link with Poisson kernel ball}), we can write
\begin{align*}
\mathbf{P}_{k,2}[f](x)&=\frac{1}{d_k}\int_{S^{d-1}}\big(P_k(x,\xi)-\mathbf{P}_{k,1}(x,\xi)\big)f(\rho\xi)\omega_k(\xi)d\sigma(\xi)\\
&=P_k[\delta_{\rho}.f](x)-\mathbf{P}_{k,1}[\delta_{\rho}.f](x).
\end{align*}
Here, $P_k[\phi]$ denotes the Poisson integral of $\phi$ and $\delta_{\rho}.f(x)=f(\rho{x})$.\\
Then, using Proposition \ref{Poisson integral part 1} and theorem A in \cite{Maslyouss}, we obtain immediately the following result:
\begin{Prop}\label{Poisson integral part 2}
Let $f$ be  a continuous function on $S(0,\rho)$.
Then $\mathbf{P}_{k,2}[f]$ is the unique solution in $\mathcal{C}^2(A)\cap\mathcal{C}(\overline{A})$ of the boundary Dirichlet problem
$$
\left\{
  \begin{array}{ll}
    \Delta_ku=0, & \hbox{on}\ A; \\
     u=0, & \hbox{on}\ S^{d-1}\\
     u=f, & \hbox{on}\ S(0,\rho).
  \end{array}
\right.
$$
\end{Prop}
\begin{Def} Let $f$ be a continuous function on $\partial A$. We define the $\Delta_k$-Poisson integral of $f$ for the annulus $A$ by
\begin{equation}\label{Poisson integral annulus}
P_{k,A}[f](x):=\mathbf{P}_{k,1}[f](x)+\mathbf{P}_{k,2}[f](x)
\end{equation}
\end{Def}
\begin{Remark}
\begin{description}
  \item[1.] Obviously  $P_{k,A}[1]=1$.
  \item[2.] Using (\ref{W-action kernel 1}) and a similar relation for the kernel $\mathbf{P}_{k,2}$, we obtain
 \begin{equation}\label{Poisson integral- W-action}
g.P_{k,A}[f]=P_{k,A}[g.f],\quad\text{with} \quad g.f(x):=f(g^{-1}x).
\end{equation}
\end{description}
\end{Remark}
From Propositions \ref{Poisson integral part 1} and \ref{Poisson integral part 2}, we deduce the following main result:

\begin{Thm}\label{Dirichlet prob annulus} Let $f\in \mathcal{C}(\partial{A})$. Then the function  $P_{k,A}[f]$ is the unique solution in $\mathcal{C}^2(A)\cap\mathcal{C}(\overline{A})$ of the boundary Dirichlet problem
$$
  \left\{
  \begin{array}{ll}
    \Delta_ku=0, & \hbox{on}\ A; \\
\\
    u=f & \hbox{on}\ \partial A.
  \end{array}
\right.
$$
\end{Thm}
From this theorem and the weak maximum principle for $\Delta_k$-harmonic function (see \cite{GalRej}), we obtain the following result:
\begin{Cor}\label{Poisson integral-harmonic function}
Let $h$ be a $\Delta_k$-harmonic function on $A$ and continuous on $\overline{A}$.
Then,
$$
\forall\ x\in A,\quad  h(x)=P_{k,A}[h](x).
$$
\end{Cor}
\section{$\Delta_k$-Green function of the annulus}
Our aim in this section is to introduce and study the Green function of the annular region $A=\{x\in\R^d,\ \rho<\|x\|<1\}$ for the Dunkl-Laplace operator. In the sequel, we will assume that $d+2\gamma>2$ i.e. $\lambda_k>0$ where $\lambda_k$ is the constant (\ref{lambda cte}).\\
 Let us first recall that  the $\Delta_k$-Newton kernel, introduced in \cite{GalRej3}, is given by
 \begin{equation}\label{Newton kernel def intr}
 N_k(x,y):=\int_0^{+\infty}p_k(t,x,y)dt,
\end{equation}
with $p_k$ the Dunkl heat kernel (see \cite{Rosler3, Rosler4})
\begin{equation}\label{heat kernel}
p_k(t,x,y)=\frac{1}{(2t)^{d/2+\gamma}c_k}\int_{\R^d}e^{-(\|x\|^2+\|y\|^2-2\prs{x}{z})/{4t}}d\mu_y(z)
\end{equation}
and $c_k$ the Macdonald-Mehta constant given by
\begin{equation*}\label{Mehta constante}
c_k:=\int_{\R^d}\exp(-\frac{\|x\|^2}{2})\omega_k(x)dx.
\end{equation*}
According to \cite{GalRej3}, the positive and symmetric kernel $N_k$ takes the following form
\begin{equation}\label{Newton kernel def}
N_k(x,y)=\frac{1}{2d_k\lambda_k}\int_{\R^d}\big(\|x\|^2+\|y\|^2-2\prs{x}{z}\big)^{-\lambda_k}d\mu_y(z).
\end{equation}
Note that if $y=0$,  then $\mu_y=\delta_0$ (with $\delta_{x_0}$ the Dirac measure at $x_0\in\R^d$) and so
$$
N_k(x,0)=\frac{1}{2d_k\lambda_k}\|x\|^{-2\lambda_k}.
$$
In addition, for each fixed $x\in \R^d$, the function $N_k(x,.)$  is $\Delta_k$-harmonic  and of class $C^{\infty}$ on  $\R^d\setminus W.x$ (where $W.x$ is the $W$-orbit of $x$),  $\Delta_k$-superharmonic (see below for precise definition) on the whole space $\R^d$ and  satisfies
$$
-\Delta_k[N_k(x,.)\omega_k]=\delta_x, \quad \text{in}\quad \mathcal{D}'(\R^d),
$$
where\\
-for  a $W$-invariant open set $\Omega\subset\R^d$, $\mathcal{D}(\Omega)$ and $\mathcal{D}'(\Omega)$ denote respectively the space of $C^{\infty}$-functions on $\Omega$ with compact support and the space of Schwartz distributions on $\Omega$.\\
-for $f\in L^1_{loc}(\Omega,\omega_k(x)dx)$, $\Delta_k(f\omega_k)$  is the Schwartz distribution on $\Omega$ defined by
$$
\prs{\Delta_k(f\omega_k)}{\varphi}=\prs{f\omega_k}{\Delta_k\varphi},\quad \varphi\in \mathcal{D}(\Omega).
$$
Moreover, for  any $x\in\R^d$, $N_k(x,x)=+\infty$.  For more details about the $\Delta_k$-Newton kernel, we refer to (\cite{GalRej3}, Section 6).\\
Let $\Omega$ be a $W$-invariant open subset of $\R^d$. Recall that  a function $u:\Omega\longrightarrow [-\infty,+\infty[$ is  $\Delta_k$-subharmonic
if (see \cite{GalRej3})
\begin{description}
  \item[1.] $u$ is upper semi-continuous on $\Omega$,
  \item[2.] $u$ is not identically $-\infty$ on each connected component of $\Omega$,
  \item[3.] $u$ satisfies the  volume sub-mean property i.e. for each closed ball  $\overline{B}(x,r)\subset\Omega$, we have
  \begin{equation}\label{volume mean def}
  u(x)\leq M_B^r(u)(x):=\frac{1}{m_k[B(0,r)]}\int_{\R^d} u(y) h_k(r,x,y)\omega_k(y) dy.
  \end{equation}
\end{description}
 Here $m_k$ is the measure $dm_k(x):=\omega_k(x)dx$ and $y\mapsto {h}_k (r,x,y)$ is the nonnegative  compactly supported measurable  function given by
 \begin{equation}\label{harmonickerneldef}
h_k(r,x,y):=\int_{\R^d}\textbf{1}_{[0,r]}(\sqrt{\|x\|^2+\|y\|^2-2\prs{x}{z}})d\mu_y(z).
\end{equation}
We refer to \cite{GalRej} for more details on the kernel $h_k$.
\bigskip

 The following result  gives some  useful facts about the Poisson integral of the $\Delta_k$-Newton kernel:
\begin{Prop}\label{properties of the Poisson-Newton} \hfill
\begin{description}
  \item[i)] For each $x\in A$, the function $P_{k,A}[N_k(x,.)]$ is the solution of the  Dirichlet problem
$$
  \left\{
  \begin{array}{ll}
    \Delta_ku=0, & \hbox{on}\ A; \\
    u=N_k(x,.) & \hbox{on}\ \partial A.
  \end{array}
\right.
$$
\item[ii)] The function $(x,y)\mapsto P_{k,A}[N_k(x,.)](y)$ is continuous on $A\times \overline{A}$.
\item[iii)] For each fixed $y\in A$, the function $x\mapsto P_{k,A}[N_k(x,.)](y)$ is $\Delta_k$-harmonic in $A$.
\end{description}
\end{Prop}
We need the following lemma:
\begin{Lem}\label{Newton kernel Continuity }
The function $(x,y)\mapsto N_k(x,y)$ is continuous on $\R^d\times\R^d\setminus\{(x,gx),\quad x\in\R^d, g\in W\}$.
\end{Lem}
\emph{Proof:} Using the following inequality (see \cite{Rosler4} Lemma 4.2)
\begin{align*}
\forall\ t>0,\quad p_k(t,x,y)\leq \frac{1}{c_k(2t)^{d/2+\gamma}}\max_{g\in W}e^{-\|x-gy\|^2/{4t}},
\end{align*}
 we can apply the Lebesgue dominated convergence theorem in formula (\ref{Newton kernel def intr}) to obtain the result of the lemma.

 \fd

\noindent\emph{Proof of Proposition \ref{properties of the Poisson-Newton}:} \textbf{i)}  If $x\in A$, then the function $N_k(x,.)$ is continuous on $\partial A$ and  by Theorem \ref{Dirichlet prob annulus}, we obtain the first assertion.
\medskip

\noindent
\textbf{ii)}
 From the first assertion, for each $x\in A$, $P_{k,A}[N_k(x,.)]$ is extendable to a continuous function on $\overline{A}$ with $P_{k,A}[N_k(x,.)]=N_k(x,.)$ on $\partial{A}$. \\
Let $(x_0,y_0)\in A\times \overline{A}$. For every $(x,y)\in A\times\overline{A}$ we have
\begin{align*}
\Big|P_{k,A}[N_k(x,.)](y)-P_{k,A}[N_k(x_0,.)](y_0)\Big|
&\leq  \Big|P_{k,A}[N_k(x,.)](y)-P_{k,A}[N_k(x_0,.)](y)\Big|\\
&+\Big|P_{k,A}[N_k(x_0,.)](y)-P_{k,A}[N_k(x_0,.)](y_0)\Big|\\
&\leq \mathbf{P}_{k,1}\Big[\big|K_{x_0}(x,.)\big|\Big](y)+ \mathbf{P}_{k,2}\Big[\big|K_{x_0}(x,.)\big|\Big](y)\\
&+\Big|P_{k,A}[N_k(x_0,.)](y)-P_{k,A}[N_k(x_0,.)](y_0)\Big|,
\end{align*}
where $K_{x_0}(x,y):=N_k(x,y)-N_k(x_0,y)$.\\
We already know that
$$
\lim_{y\rightarrow y_0}P_{k,A}[N_k(x_0,.)](y)=P_{k,A}[N_k(x_0,.)](y_0).
$$
 Now, let $\varepsilon>0$ and $R>0$ be such that $\overline{B}(x_0,R)\subset{A}$.  Since $(x,\xi)\longmapsto N_k(x,\xi)$ is uniformly continuous on $\overline{B}(x_0,R)\times S^{d-1}$, we deduce that there exists $\eta>0$ such that
   $$
   \forall\ (x,\xi)\in B(x_0,\eta)\times S^{d-1},\quad |K_{x_0}(x,\xi)|=\big|N_k(x,\xi)-N_k(x_0,\xi)\big|<\varepsilon.
   $$
 Then, using (\ref{solution1}) as well as the inequalities (\ref{Poisson kernel as prob measure}) and (\ref{norme1 Poisson kernel annulus}), we get for every $x \in B(x_0,\eta)$ and every $y\in A$
 $$
 \mathbf{P}_{k,1}\big[\big|K_{x_0}(x,.)\big|\big](y)\leq \frac{1}{d_k}\int_{S^{d-1}}\mathbf{P}_{k,1}(y,\xi)|K_{x_0}(x,\xi)|\omega_k(\xi)d\sigma(\xi)\leq \varepsilon.
 $$
The same idea works if we replace the kernel $\mathbf{P}_{k,1}$ by  $\mathbf{P}_{k,2}$. Finally, we obtain
$$
\lim_{(x,y)\rightarrow (x_0,y_0)}P_{k,A}[N_k(x,.)](y)=P_{k,A}[N_k(x_0,.)](y_0).
$$
That is the function $(x,y)\mapsto P_{k,A}[N_k(x,.)](y)$ is continuous on $A\times \overline{A}$ as desired.
\medskip

\noindent
\textbf{iii)} According to Corollary 4.6 in \cite{GalRej2}, it is enough to show that the functions $x\mapsto u_y(x):=\textbf{P}_{k,1}[N_k(x,.)](y)$ and $x\mapsto v_y(x):=\textbf{P}_{k,2}[N_k(x,.)](y)$ satisfy  the volume-mean property.\\
 Let then $x_0\in A$ and $R>0$ such that $\overline{B}(x_0,R)\subset A$. As the kernels $N_k$, $h_k$ \emph{}and $\textbf{P}_{k,1}$ are nonnegative, we can use Fubini's theorem to obtain
\begin{align*}
M_B^R(u_y)(x_0)=\frac{1}{d_k}\int_{S^{d-1}}\textbf{P}_{k,1}(y,\xi)M_B^R[N_k(.,\xi)](x_0)\omega_k(\xi)d\sigma(\xi)
\end{align*}
But for any $\xi\in S^{d-1}$, the function $N_k(.,\xi)$ is $\Delta_k$-harmonic on $A$. Hence, it satisfies the volume-mean property i.e.  $M_B^R[N_k(.,\xi)](x_0)=N_k(x_0,\xi)$. Therefore, we obtain
\begin{align*}
M_B^R(u_y)(x_0)&=\frac{1}{d_k}\int_{S^{d-1}}\textbf{P}_{k,1}(y,\xi)N_k(x_0,\xi)
\omega_k(\xi)d\sigma(\xi)\\
&=\textbf{P}_{k,1}[N_k(x_0,.)](y)=u_y(x_0).
\end{align*}
We prove similarly that   $x\mapsto v_y(x):=\textbf{P}_{k,2}[N_k(x,.)](y)$ is also a $\Delta_k$-harmonic function in $A$.
This shows the desired result.

\fd

\begin{Def} For $x\in A$, the function $G_{k,A}(x,.)$ defined by
\begin{equation}\label{Green function def}
G_{k,A}(x,y):=N_k(x,y)-P_{k,A}[N_k(x,.)](y),\quad y\in A,
\end{equation}
is called the $\Delta_k$-Green function of $A$ with pole $x$.
\end{Def}
The $\Delta_k$-Green function $G_{k,A}$ has the following properties:
\begin{Prop} \label{Green function properties} Let $x\in A$. Then \hfill
\begin{description}
\item[1.] The function $G_{k,A}(x,.)$ is $\Delta_k$-harmonic on $A\setminus W.x$, is $\Delta_k$-superharmonic on $A$ and satisfies
\begin{equation}\label{Green-Poisson equation1}
-\Delta_k[G_{k,A}(x,.)\omega_k]=\delta_x\quad \text{in}\quad \mathcal{D}'(A).
\end{equation}
\item[2.] $G_{k,A}(x,x)=+\infty$ and $G_{k,A}(x,y)<+\infty$ whenever $y\notin W.x$ .
\item[3.]  For every $\xi\in \partial{A}$, $\lim_{y\rightarrow\xi}G_{k,A}(x,y)=0$.
\item[4.] For every $y\in A$, $G_{k,A}(x,y)>0$.
\item[5.] For every $x,y\in A$, $G_{k,A}(x,y)=G_{k,A}(y,x)$.
\item[6.] For every $x,y\in A$ and $g\in W$, $G_{k,A}(gx,gy)=G_{k,A}(x,y)$.
\item[7.] The zero function  is the greatest $\Delta_k$-subharmonic minorant of  $G_{k,A}(x,.)$ on $A$.
\item[8.] The function $(x,y)\mapsto G_{k,A}(x,y)$ is continuous on $A\times \overline{A}\setminus\{(x,gx):\ x\in A,g\in W\}$.
\end{description}
\end{Prop}
\noindent \emph{ Proof: }
The first and the second assertions follow from the properties of the $\Delta_k$-Newton kernel previously mentioned. In addition, by Proposition \ref{properties of the Poisson-Newton}, we easily obtain the third statement.
\medskip

\noindent
\textbf{4)} As $\lim_{y\rightarrow\xi\in\partial{A}}G_{k,A}(x,y)=0$, the weak minimum principle for $\Delta_k$-superharmonic functions (see \cite{GalRej3}, Theorem 3.1) implies that $G_{k,A}(x,.)\geq 0$ on $A$.\\
If $G_{k,A}(x,y_0)=0$ for some $y_0\in A$, it follows from the strong maximum principle (see \cite{GalRej3}) that $G_{k,A}(x,.)$ is the zero function on $A$ which is impossible because $G_{k,A}(x,x)=+\infty$. Thus, $G_{k,A}$ is a positive kernel on $A\times {A}$.
\medskip

\noindent
\textbf{5)} Since $N_k$ is a symmetric  kernel, we have to prove that
$$
\forall\ x,y\in A, \quad P_{k,A}[N_k(x,.)](y)= P_{k,A}[N_k(y,.)](x).
$$
 Let $y\in A$ and consider the function
\begin{align*}
H_y(x):=P_{k,A}[N_k(x,.)](y)-P_{k,A}[N_k(y,.)](x).
\end{align*}
From Proposition \ref{properties of the Poisson-Newton} , i) and iii), $H_y$ is a $\Delta_k$-harmonic function in $A$. On the other hand, writing
$$
H_y(x)=N_k(x,y)-G_{k,A}(x,y)-P_{k,A}[N_k(y,.)](x)
$$
and using the positivity of $G_{k,A}$ as well as the symmetry property of the kernel $N_k$, we conclude that
$$
\limsup_{x\mapsto \xi\in\partial{A}}H_y(x)\leq N_k(\xi,y)-N_k(y,\xi)=0.
$$
Then the weak maximum principle yields that $H_y\leq 0$ on $A$. That is, we have
$$
\forall\ x,y\in A, \quad P_{k,A}[N_k(x,.)](y)\leq P_{k,A}[N_k(y,.)](x).
$$
By interchanging the role of $x$ and $y$, we also get the reverse inequality.  Finally, we obtain the desired equality.
\medskip

\noindent
\textbf{6)} The result follows immediately from (\ref{Poisson integral- W-action}) and from the relation  $N_k(gx,gy)=N_k(x,y)$, $x,y\in \R^d$, $g\in W$ (see \cite{GalRej3}).
\medskip

\noindent
\textbf{7)} As $G_{k,A}(x,.)$ is positive on $A$, we know that the zero function  is a $\Delta_k$-subharmonic minorant of $G_{k,A}(x,.)$.\\
 Now, let $s$ be a $\Delta_k$-subharmonic function on $A$ such that $s\leq G_{k,A}(x,.)$ on $A$.
Using  the statement 3), we obtain $\limsup_{z\rightarrow\xi\in\partial A}s(z)\leq0$.
Thus the weak maximum principle for $\Delta_k$-subharmonic functions (see \cite{GalRej3})  yields that
$s\leq 0$ on $A$.
\medskip

\noindent
\textbf{8)} The result follows immediately from the statement ii)  of Proposition \ref{properties of the Poisson-Newton} and Lemma \ref{Newton kernel Continuity }.
\fd

\medskip

In the following result, we will express the Green function $G_{k,A}$ in terms of the $\Delta_k$-spherical harmonics. More  precisely, we have
\begin{Thm}\label{Green function-zonal thm} The $\Delta_k$-Green function in $A$ is given by
\begin{equation}\label{Green function-zonal}
G_{k,A}(x,y)=N_k(x,y)-\sum_{n=0}^{+\infty}\frac{a_{k,n}(y)\|x\|^n+b_{k,n}(y)\|x\|^{-n-2\lambda_k}\rho^n}{d_k(2\lambda_k+2n)}Z_{k,n}\big(\frac{x}{\|x\|},y\big).
\end{equation}
\end{Thm}
We need the following result:
\begin{Prop} For $x,y\in\R^d$ such that $\|y\|<\|x\|$, we have
\begin{equation}\label{Newton kenel Gegenbauer}
N_k(x,y)=\sum_{n=0}^{+\infty}\frac{\|x\|^{-2\lambda_k}}{d_k(2\lambda_k+2n)}\|y\|^n\|x\|^{-n}Z_{k,n}\big(\frac{x}{\|x\|},\frac{y}{\|y\|}\big).
\end{equation}
\end{Prop}
\emph{Proof:} Let $\|y\|<\|x\|$. From (\ref{Newton kernel def}), we have
\begin{align*}
N_k(x,y)&=\frac{\|x\|^{-2\lambda_k}}{2d_k\lambda_k}\int_{\R^d}\big(1-\frac{2\prs{x}{z}}{\|x\|^2}+
\frac{\|y\|^2}{\|x\|^2}\big)^{-\lambda_k}d\mu_y(z)\\
&=\frac{\|x\|^{-2\lambda_k}}{2d_k\lambda_k}\int_{\R^d}\sum_{n=0}^{+\infty}\|y\|^n\|x\|^{-n}\frac{(2\lambda_k)_n}{n!}P_n^{\lambda_k}
\Big(\prs{\frac{x}{\|x\|}}{\frac{z}{\|y\|}}\Big)d\mu_y(z)\\
&=\sum_{n=0}^{+\infty}\frac{\|x\|^{-2\lambda_k}}{d_k(2\lambda_k+2n)}\|y\|^n\|x\|^{-n}Z_{k,n}\big(\frac{x}{\|x\|},\frac{y}{\|y\|}\big);
\end{align*}
where in the second line, we have used the relation (\ref{support mu x}) and the generating relation  (see for example \cite{DunklXu}, p. 18)
$$
(1-2ar+r^2)^{-\mu}=\sum_{n=0}^{+\infty}\frac{(2\mu)_n}{n!}P_n^{\mu}(a)r^n,\quad\mu>0,\quad |r|<1,\ |a|\leq1
$$
and in the last line, we have used \\
- the inequality $\underset{x\in[-1,1]}{\sup}|P_n^{\mu}(x)|\leq P_n^{\mu}(1)$ (see \cite{Szego}, Theorem 7.32.1) and the above generation relation with $a=1$  which allows us to permute the symbols $\sum$ and $\int$,\\
- the fact that $\mu_{\frac{y}{\|y\|}}$ is the image measure of $\mu_y$ by the dilation $\xi\mapsto \frac{\xi}{\|y\|}$\\
- the relations (\ref{Integral rep Vk }) and (\ref{sphercial harrmonics Gegen}).
\fd

\medskip

\noindent\emph{Proof of Theorem \ref{Green function-zonal thm}:}  By Theorem \ref{Dirichlet prob annulus}, we have
\begin{align*}
P_{k,A}[N_k(x,.)](y)=\mathbf{P}_{k,1}[N_k(x,.)](y)+ \mathbf{P}_{k,2}[N_k(x,.)](y):=I_1+I_2.
\end{align*}
$\bullet$ We have
\begin{align*}
I_1&=\sum_{n=0}^{+\infty}\frac{a_{k,n}(y)}{d_k}\int_{S^{d-1}}Z_{k,n}(y,\xi)N_k(x,\xi)\omega_k(\xi)d\sigma(\xi)\\
&=\sum_{n=0}^{+\infty}\frac{a_{k,n}(y)}{d_k}\int_{S^{d-1}}\sum_{m=0}^{+\infty}\frac{\|x\|^m}{d_k(2\lambda_k+2m)}
Z_{k,m}\big(\frac{x}{\|x\|},\xi\big)Z_{k,n}(y,\xi)\omega_k(\xi)d\sigma(\xi)\\
&=\sum_{n=0}^{+\infty}a_{k,n}(y)\sum_{m=0}^{+\infty}\frac{\|x\|^m}{d_k(2\lambda_k+2m)}\frac{1}{d_k}\int_{S^{d-1}}
Z_{k,m}\big(\frac{x}{\|x\|},\xi\big)Z_{k,n}(y,\xi)\omega_k(\xi)d\sigma(\xi)\\
&=\sum_{n=0}^{+\infty}\frac{a_{k,n}(y)\|x\|^n}{d_k(2\lambda_k+2n)}Z_{k,n}\big(\frac{x}{\|x\|},y\big),
\end{align*}
where, we have used \\
-the relation (\ref{Newton kenel Gegenbauer}) in the second line;\\
-the inequalities (\ref{k-sphercial harmonic ineq}) and  $\|x\|<1$ in order to interchange the symbols  $\int$ and $\sum$ in the third line;\\
-the fact that  $\mathcal{H}_{\Delta_k,n}(\R^d)\perp\mathcal{H}_{\Delta_k,m}(\R^d)$ whenever $m\neq n$ and the reproducing formula (\ref{reproducing kernel formula}) in the last line.
\medskip

\noindent Notice that if $\|x\|>\rho$, then (\ref{Newton kenel Gegenbauer}) yields that
$$
\forall\ \xi \in S^{d-1},\quad N_k(x,\rho\xi)=\sum_{n=0}^{+\infty}\frac{\|x\|^{-2\lambda_k}}{d_k(2\lambda_k+2n)}\rho^n\|x\|^{-n}Z_{k,n}\big(\frac{x}{\|x\|},\xi\big).
$$
Hence we obtain similarly
\begin{align*}
I_2&=\sum_{n=0}^{+\infty}\frac{b_{k,n}(y)}{d_k}\int_{S^{d-1}}Z_{k,n}(y,\xi)N_k(x,\rho\xi)\omega_k(\xi)d\sigma(\xi)\\
&=\sum_{n=0}^{+\infty}\frac{b_{k,n}(y)}{d_k}\int_{S^{d-1}}\sum_{m=0}^{+\infty}\frac{\|x\|^{-2\lambda_k-m}\rho^m}{d_k(2\lambda_k+2m)}
Z_{k,m}\big(\frac{x}{\|x\|},\xi\big)Z_{k,n}(y,\xi)\omega_k(\xi)d\sigma(\xi)\\
&=\sum_{n=0}^{+\infty}b_{k,n}(y)\sum_{m=0}^{+\infty}\frac{\|x\|^{-2\lambda_k-m}\rho^m}{d_k(2\lambda_k+2m)}\frac{1}{d_k}\int_{S^{d-1}}
Z_{k,m}\big(\frac{x}{\|x\|},\xi\big)Z_{k,n}(y,\xi)\omega_k(\xi)d\sigma(\xi)\\
&=\sum_{n=0}^{+\infty}\frac{b_{k,n}(y)\|x\|^{-2\lambda_k-n}\rho^n}{d_k(2\lambda_k+2n)}Z_{k,n}\big(\frac{x}{\|x\|},y\big).
\end{align*}
This gives the desired formula (\ref{Green function-zonal}).
\fd
\begin{Remark} Using (\ref{Newton kenel Gegenbauer}) and replacing the functions $a_{k,n}$ and $b_{k,n}$ by their expressions, if $x,y\in A$ with $\|y\|<\|x\|$ we can write
\begin{align*}
G_{k,A}(x,y)=\sum_{n=0}^{+\infty}\frac{\Big(\|y\|^{2\lambda_k+2n}-\rho^{2\lambda_k+2n}\Big)(1-\|x\|^{2\lambda_k+2n}\Big)}
{d_k(2\lambda_k+2n)(1-\rho^{2\lambda_k+2n})(\|x\|\|y\|)^{2\lambda_k+n}} Z_{k,n}\big(\frac{x}{\|x\|},\frac{y}{\|y\|}\big).\\
\end{align*}
This formula generalizes the classical case (if $k=0$, $2\lambda_0=d-2$) proved in \cite{Grossi}.
\end{Remark}
\section{Applications}
\subsection{Poisson-Jensen formula for  $\Delta_k$-subharmonic functions}
Our goal now is to prove an analogue of the Poisson-Jensen formula for  $\Delta_k$-subharmonic functions on $\Omega\supset \overline{A}$. Note that a Poisson-Jensen formula has been proved when $u$ is a $C^2-\Delta_k$-subharmonic function on $\Omega$ which  contains the closed unit ball (see \cite{Graczyk-Luks-Rosler}).
\begin{Thm}\label{Poisson jensen formula} Let $u$ be a $\Delta_k$-subharmonic function on a $W$-invariant open set $\Omega\supset\overline{A}$. Then,
\begin{equation}
u(x)=P_{k,A}[u](x)-\int_{A}G_{k,A}(x,y)d\nu_u(y),\quad x\in A,
\end{equation}
where $\nu_u:=\Delta_k(u\omega_k)$ is the $\Delta_k$-Riesz measure of $u$ (see \cite{GalRej3}).
\end{Thm}
\emph{Proof:} Let $O$ be a bounded $W$-invariant open set such that $\overline{A}\subset O\subset \overline{O}\subset \Omega$ . Using the Riesz decomposition theorem for $\Delta_k$-subharmonic functions (see \cite{GalRej3}), we deduce that there exists a $\Delta_k$-harmonic function $h$ on $O$ such that
$$
\forall\ x\in O,\quad u(x)=h(x)-\int_{O}N_k(x,y)d\nu_u(y):=h(x)-s(x).
$$
Then, we have
$$
\forall\ x\in A, \quad P_{k,A}[u](x)=P_{k,A}[h](x)-P_{k,A}[s](x).
$$
From Corollary \ref{Poisson integral-harmonic function}, we have $P_{k,A}[h]=h$ on $A$. Moreover, for $x\in A$, we have
$$
P_{k,A}[s](x)=\mathbf{P}_{k,1}[s](x)+\mathbf{P}_{k,2}[s](x).
$$
The crucial part here is to show that
\begin{equation}\label{crucial part Poisson-Jensen}
\forall\ x\notin A,\quad P_{k,A}[N_k(x,.)]=N_k(x,.)\quad \text{on}\ A.
\end{equation}
Assume this relation for the moment. By Fubini's theorem, we have
\begin{align*}
\mathbf{P}_{k,1}[s](x)&:=\frac{1}{d_k}\int_{S^{d-1}} \mathbf{P}_{k,1}(x,\xi)s(\xi)\omega_k(\xi)d\sigma(\xi)\\
&=\frac{1}{d_k}\int_{O}\int_{S^{d-1}} \mathbf{P}_{k,1}(x,\xi)N_k(\xi,y)\omega_k(\xi)d\sigma(\xi)d\nu_u(y)\\
&=\int_{O}\mathbf{P}_{k,1}[N_k(y,.)](x)d\nu_u(y).
\end{align*}
By the same way, we also have
$$
\mathbf{P}_{k,2}[s](x)=\int_{O}\mathbf{P}_{k,2}[N_k(y,.)](x)d\nu_u(y).
$$
The above relations as well as (\ref{crucial part Poisson-Jensen}) imply that
\begin{align*}
P_{k,A}[s](x)&=\int_{O}P_{k,A}[N_k(y,.)](x)d\nu_u(y)\\
&=\int_{A}P_{k,A}[N_k(y,.)](x)d\nu_u(y)+\int_{O\setminus A}P_{k,A}[N_k(y,.)](x)d\nu_u(y)\\
&=\int_{A}\big(N_k(x,y)-G_{k,A}(x,y)\big)d\nu_u(y)+\int_{O\setminus A}N_k(x,y)d\nu_u(y)\\
&=\int_{O}N_k(x,y)d\nu_u(y)-\int_{A}G_{k,A}(x,y)d\nu_u(y).
\end{align*}
This implies the desired Poisson-Jensen formula. Now, it remains to prove (\ref{crucial part Poisson-Jensen}).  We will distinguish three cases.
\medskip

\noindent \emph{First case:} $x\notin \overline{A}$. As $N_k(x,.)$ is $\Delta_k$-harmonic on $A$ and continuous on $\overline{A}$, we deduce by Corollary \ref{Poisson integral-harmonic function} that $P_{k,A}[N_k(x,.)]=N_k(x,.)$ on $A$.
\medskip

\noindent\emph{Second case:} $x\in S^{d-1}$. For $\varepsilon>0$ small enough, the function $N_k\big((1+\varepsilon)x,.\big)$ is $\Delta_k$-harmonic in the open ball  $B(0,1+\varepsilon)\supset \overline{A}$.
Therefore, again by Corollary \ref{Poisson integral-harmonic function}, we obtain
\begin{equation}\label{case 2}
\forall\ y\in A,\quad N_k\big((1+\varepsilon)x,y\big)=P_{k,A}[N_k\big((1+\varepsilon)x,.\big)](y).
\end{equation}
  Clearly we have $\lim_{\varepsilon\rightarrow0} N_k\big((1+\varepsilon)x,y\big)=N_k(x,y)$ for every fixed $y\in A$. Moreover, using (\ref{Newton kernel def}) and the fact that supp$\ \mu_y\subset \overline{B}(0,\|y\|)$ we can see that
$$
N_k\big((1+\varepsilon)x,y\big)\leq N_k(x,y),\quad\text{whenever}\quad \|y\|\leq \|x\|.
$$
Consequently, we can use  the Lebesgue dominated convergence theorem to obtain
$$
\forall\ y\in A,\quad \lim_{\varepsilon\rightarrow0}P_{k,A}[N_k\big((1+\varepsilon)x,.\big)](y)=P_{k,A}[N_k(x,.)](y).
$$
Hence, letting $\varepsilon\longrightarrow 0$ in the relation (\ref{case 2}), we get the result in this case.
\medskip

\noindent
 \emph{Third case:}  $x\in S(0,\rho)$. Let $0<\varepsilon<1/2$. In this case, the function $N_k\big((1-\varepsilon)x,.\big)$ is $\Delta_k$-harmonic in
 $\R^d\setminus \overline{B}\big(0,(1-\varepsilon)\rho\big)\supset \overline{A}$ and then from Corollary \ref{Poisson integral-harmonic function} we deduce that
$$
\forall\ y\in A,\quad N_k\big((1-\varepsilon)x,y\big)=P_{k,A}[N_k\big((1-\varepsilon)x,.\big)](y).
$$
Note that from (\ref{support mu x}), we can write
$$
N_k(x,y)=\frac{1}{2d_k\lambda_k}\int_{\R^d}\Big(\sum_{g\in W}\lambda_g(z)\|x-gy\|^2\Big)^{-\lambda_k}d\mu_y(z),
$$
where for every $z\in \text{supp}\ \mu_y$, the nonnegative numbers $\lambda_g(z)$ are such that $\sum_{g\in W}\lambda_g(z)=1$.
Using the above relation we easily see that
$$
N_k\big((1-\varepsilon)x,y\big)\leq 2^{2\lambda_k}N_k(x,y),\quad\text{whenever}\quad \|x\|\leq \|y\|,
$$
Finally by same way as in the second case we obtain the result.

\fd
\subsection{Positive solution of $\Delta_k$-nonlinear elliptic  problem on the annulus}
In this section, we will investigate the positive continuous solutions of the semilinear problem
\begin{equation*}
\Delta_k(u\omega_k)=\phi(.,u)\omega_k\quad \text{in}\quad \mathcal{D}'(A),
\end{equation*}
in the sense that
$$
\forall\ \varphi\in \mathcal{D}(A),\quad \int_A u(x)\Delta_k\varphi(x)\omega_k(x)dx=\int_A\varphi(x)\phi(x,u(x))\omega_k(x)dx.
$$
We  assume  that  $\phi(x,u(x))=\phi_1(x)\phi_2(u(x))$ where\\
 $\bullet$ $\phi_1$ is a  nonnegative bounded measurable function on $A$.\\
 $\bullet$ $\phi_2$ is a  nonnegative and nondecreasing continuous function on $[0,+\infty[$ with $\phi_2(0)=0$.\\

 In \cite{Ben Chrouda 2}, by using some tools from probabilistic potential theory,  the authors have studied the positive solution on the unit ball  $B$  of the semilinear problem
  $$
  \Delta_k(u)=\varphi(u)\quad \text{in}\ \mathcal{D}'(B)\quad \text{and}\quad u=f\quad \text{on}\ \partial{B}.
  $$
Let us denote by $\mathcal{C}^+(\overline{A})$ the convex cone of nonnegative and continuous functions on $A$.
\begin{Thm}\label{thm semilinear problem}
	Let $\phi$ be as above. Then, for every  $f\in\mathcal{C}^+( \partial{A})$, the semilinear Dirichlet problem
\begin{equation}\label{semilinear-Dirichlet_problem}
	\left\{
	\begin{array}{ll}
	\Delta_k(u\omega_k)=\phi(.,u)\omega_k, &\mbox{ in }\ \mathcal{D}'(A)\\
\\
	u= f, &\mbox{ on }\  \partial A
\end{array}
	\right.
	\end{equation}
admits one and only one solution $u \in \mathcal{C}^+(\overline{A})$. Furthermore, we have
$$
 \forall\ x\in A,\quad u(x) + \int_AG_{k,A}(x,y)\phi(y,u(y))\omega_k(y)dy= P_{k,A}[f](x).
$$	
\end{Thm}	
We begin by showing the uniqueness of the solution. This fact follows immediately from the following maximum principle type result:
 \begin{Lem}\label{Dunkl comparaison principle}
Let $u,v\in \mathcal{C}(A)$ and let $\phi$ be a function satisfying the above conditions. If
\begin{equation*}
\left\{
\begin{array}{ll}
\Delta_k (u\omega_k)- \phi(.,u)\omega_k\leq\Delta_k(v\omega_k) - \phi(.,v)\omega_k , & \hbox{in}\  \mathcal{D}'(A), \\
\\
\limsup\limits_{x \to y\in\partial{A}} (v-u) (x) \leq 0, &
\end{array}
\right.
\end{equation*}
then $v\leq u$ in $A$.
\end{Lem}	
\emph{Proof:} Let  $U$ be the upper semi-continuous function  defined  by
 $$
U(x)=\left\{
       \begin{array}{ll}
         v(x)-u(x), & \hbox{if}\ x\in A; \\
         \limsup\limits_{y \to x\in\partial{A}} (v-u)(y), & \hbox{if}\ x\in \partial{A}
       \end{array}
     \right.
$$
and $x_0\in\overline{A}$ be such that $U(x_0)=\max_{\overline{A}}U$.\\
We have to prove that $U(x_0)\leq0$. We suppose the contrary i.e.  $U(x_0)>0$. As $U\leq0$ on $\partial{A}$, this implies that $x_0\notin \partial{A}$.\\
Let $O$ be the nonempty open set given by
 $$
 O:=\{x\in A:\quad U(x)>0\}
 $$
and $\Omega$ be the connected component of $x_0$ in $O$ which is also an open set of $\R^d$. \\
To get a contradiction, we claim that it is suffices  to establish  that
\begin{equation}\label{main idea lemma unicite}
U=U(x_0)\quad \text{on}\quad \Omega.
\end{equation}
Indeed,\\
$\bullet$ If $\Omega=O$ (i.e. $O$ is connected), then (\ref{main idea lemma unicite}) holds on $O=\Omega$.
But
$$
\partial{O}\subset \partial{A}\cup (A\setminus{O})=\{x\in \overline{A},\ U(x)\leq0\}.
$$
Consequently, using the fact that $U$ is upper semi-continuous and (\ref{main idea lemma unicite}),  we get
 $$
 \forall\ x\in \partial{O},\quad U(x)=\limsup_{y\rightarrow{x},y\in {O}}U(y)=U(x_0).
 $$
 Thus, we obtain a contradiction. \\
$\bullet$ If $\Omega\neq O$, then as $\Omega$ is a connected component of $O$ we have   $\partial{\Omega}\cap O=\emptyset$. Therefore,  we have $\partial{\Omega}\subset\partial{A}\cup A\setminus{O}$ and as above we get a contradiction.
\medskip

Now, our aim is  to  prove that $U=U(x_0)$ on $\Omega$. For this, we introduce the nonempty closed set
$$
\Omega_0:=\{x\in \Omega:\quad U(x)=U(x_0)\}.
$$
Note that
$$
\Omega=\big(\Omega\cap \cup_{\alpha\in \Rac} H_{\alpha}\big)\cup \big(\Omega\setminus \cup_{\alpha\in \Rac} H_{\alpha}\big)=\big(\Omega\cap\cup_{\alpha\in \Rac} H_{\alpha}\big)\cup\big(\cup_{g\in W}\Omega \cap{g}.\mathbf{C}\big),
$$
where $\mathbf{C}$ is a fixed Weyl chamber and  $g\mathbf{C}$, $g\in W$, are the connected components of $\R^d\setminus \cup_{\alpha\in \Rac}H_{\alpha}$.\\
Fix $\xi\in \Omega_0$ and $R>0$ such that the open ball $B(\xi,R)$ is contained in $\Omega$.  We will distinguish three  possible locations of $\xi$ depending on the sets $E_1:=\{\alpha\in\Rac,\quad U(\sigma_{\alpha}\xi)=U(\xi)\}$ and $E_2:=\{\alpha\in\Rac,\quad \xi\in H_{\alpha}\}\subset E_1$.
\medskip

\noindent \textbf{First case:} $E_1=\Rac$. This implies that $U(g\xi)=U(\xi)>0$ for all $g\in W$. Moreover, clearly there exists $r\in (0,R]$ such that
 \begin{equation}\label{first case}
 \forall\ g\in W,\ \forall\ x\in B(g\xi,r),\quad U(x)\geq0.
 \end{equation}
 Consider the $W$-invariant continuous function $U^W$ defined on $A$ by
$$
U^W(x):=\frac{1}{|W|}\sum_{g\in W}g.U(x)=\frac{1}{|W|}\sum_{g\in W}U(g^{-1}x).
$$
 We easily see that $U^W$ has a maximum at the point $\xi$ with $U^W(\xi)=U(\xi)=U(x_0)$. Furthermore, using the $W$-invariance property of $\Delta_k$ (i.e. $g\circ\Delta_k=\Delta_k\circ g$) and the hypothesis of the lemma, we obtain
\begin{align*}
\Delta_k(U^W\omega_k)&=\frac{1}{|W|}\sum_{g\in W}g.[\Delta_k(U\omega_k)]\geq \frac{1}{|W|}\sum_{g\in W}g.[\big(\phi(.,v)-\phi(.,u)\big)\omega_k]\quad \text{in}\quad \mathcal{D}'(A).
\end{align*}
Now, since $U\geq0$ on $B^W(\xi,r)=\cup_{g\in W}B(g\xi,r)$ (from (\ref{first case})) and  $\phi_2$ is nondecreasing, we deduce that
$$
\Delta_k(U^W\omega_k)\geq 0\quad \text{in}\quad \mathcal{D}'(B^W(\xi,r)).
$$
That is $U$ is weakly $\Delta_k$-subharmonic on $B^W(\xi,r)$. But the continuity of $U$ and the Weyl lemma for $\Delta_k$-subharmonic functions (see \cite{GalRej3}, Theorem 5.2) imply that  $U^W$ is strongly $\Delta_k$-subharmonic on the open $W$-invariant set $B^W(\xi,r)$.\\
Now, if we follow  the proof of the strong maximum principle in \cite{GalRej3} for  the $W$-invariant $\Delta_k$-subharmonic function $U^W$, then we conclude that
$$
 U=U(\xi)=U(x_0),\quad \text{on}\quad B(\xi,r).
$$
Hence, we have $B(\xi,r)\subset \Omega_0$.
\medskip

\noindent
\textbf{Second case:}\ $E_1\neq\Rac$ and $E_2=\emptyset$ i.e. $\xi\notin\cup_{\alpha\in \Rac} H_{\alpha}$. So there is a unique $g_0\in W$ such that  $\xi\in \Omega\cap g_0\mathbf{C}$. Clearly, we can suppose that $B(\xi,R)\subset \Omega\cap g_0\mathbf{C}$.\\
Let $U^W$ be the W-invariant continuous function defined on $B^W(\xi,R):=\cup_{g\in W}B(g\xi,R)$ by
$$
U^W(x)=g.U(x):=U({g}^{-1}.x)\quad \text{whenever}\quad x\in B(g\xi,R).
$$
We are going to establish that the  function $U^W$ is $\Delta_k$-subharmonic on $B^W(\xi,r)$ for some $r>0$ will be chosen later. Again from the continuity of
$U^W$  and the Weyl lemma, it is enough  to show that $U^W$ is $\Delta_k$-subharmonic in sense of distributions.\\
$\bullet$ Firstly, we have the following decomposition
 $$
 U^W=\sum_{i=1}^nU^W\textbf{1}_{B(g_i\xi,R)}=\sum_{i=1}^ng_i.[U\textbf{1}_{B(\xi,R)}],
$$
where $g_1=id,g_2,\dots,g_n\in W$ are such that
$
\textbf{1}_{B^W(\xi,R)}=\sum_{i=1}^n\textbf{1}_{B(g_i\xi,R)}$.\\
$\bullet$ Secondly, for $f\in \mathcal{C}^2(B^W(\xi,R))$, we  can write  $\Delta_k=L_k-A_k$ where
 \begin{equation*}
L_kf(x)= \Delta f(x)+2\sum_{\alpha\in \Rac_+}k(\alpha)\frac{\prs{\nabla f(x)}{\alpha}}{\prs{\alpha}{x}}
\end{equation*}
and
$$
A_kf(x)=\sum_{\alpha\in \Rac_+}k(\alpha)\|\alpha\|^2\frac{f(x)-f(\sigma_{\alpha}(x))}{\prs{\alpha}{x}^2}.
$$
$\bullet$ For $\varphi\in\mathcal{D}(B^W(\xi,R))$ nonnegative, we have
\begin{align*}
\prs{\Delta_k(U^W\omega_k)}{\varphi}&=\sum_{i=1}^n\prs{\Delta_k\big(g_i.[U\textbf{1}_{B(\xi,R)}]\omega_k\big)}{\varphi}
=\sum_{i=1}^n\prs{g_i.[U\omega_k\textbf{1}_{B(\xi,R)}]}{\Delta_k\varphi}\\
&=\sum_{i=1}^n\prs{g_i.[U\omega_k\textbf{1}_{B(\xi,R)}]}{L_k\varphi-A_k\varphi}\\
&=\sum_{i=1}^n\prs{U\omega_k}{(L_k[g_i^{-1}.\varphi])\textbf{1}_{B(\xi,R)}}-\prs{\sum_{i=1}^ng_i.[U\textbf{1}_{B(\xi,R)}]\omega_k}{A_k(\varphi)}\\
&=\sum_{i=1}^n\prs{U\omega_k}{L_k\big([g_i^{-1}.\varphi]\textbf{1}_{B(\xi,R)}\big)}-\underbrace{\prs{U^W\omega_k}{A_k(\varphi)}}_{=0}\\
&=\sum_{i=1}^n\prs{U\omega_k}{\Delta_k\big([g_i^{-1}.\varphi]\textbf{1}_{B(\xi,R)}\big)+A_k([g_i^{-1}.\varphi]\textbf{1}_{B(\xi,R)}\big)}\\
&=\sum_{i=1}^n\prs{\Delta_k(U\omega_k)}{[g_i^{-1}.\varphi]\textbf{1}_{B(\xi,R)}}
+\sum_{i=1}^n\prs{U\omega_k}{A_k\big([g_i^{-1}.\varphi]\textbf{1}_{B(\xi,R)}\big)}\\
&\geq \sum_{i=1}^n\prs{[\phi(.,v)-\phi(.,u)]\omega_k}{[g_i^{-1}.\varphi]\textbf{1}_{B(\xi,R)}}
+\sum_{i=1}^n\prs{U\omega_k}{A_k\big([g_i^{-1}.\varphi]\textbf{1}_{B(\xi,R)}\big)}\\
&\geq \sum_{i=1}^n\prs{U\omega_k}{A_k\big([g_i^{-1}.\varphi]\textbf{1}_{B(\xi,R)}\big)},
\end{align*}
where we have used \\
- the fact that $L_k$ commutes with the $W$-action i.e. $L_k\circ g=g\circ L_k$, $g\in W$, in the third line
and the fact that it  preserves the support in the forth line,\\
-  the $W$-invariance property of $U^W$ which implies that $\prs{U^W\omega_k}{A_k(\varphi)}=0$ in the forth line,\\
-  the decomposition $L_k=\Delta_k+A_k$ in the fifth line,\\
- the fact that $[g_i^{-1}.\varphi]\textbf{1}_{B(\xi,R)}\in \mathcal{D}(B(\xi,R))$ in the sixth line,\\
- the hypothesis of the lemma in the seventh line,\\
- the nondecreasing property of $\phi_2$ in the last line. \\
$\bullet$ As $A_k$ is a symmetric operator in the sense that $\prs{A_k(f)\omega_k}{\psi}=\prs{f\omega_k}{A_k(\psi)}$, it yields that
 $$
 \prs{\Delta_k(U^W\omega_k)}{\varphi}\geq \sum_{i=1}^n\prs{A_k(U)\omega_k}{[g_i^{-1}.\varphi]\textbf{1}_{B(\xi,R)}},
 $$
 Clearly $A_k(U)(\xi)\geq0$. But, since $E_1\neq \Rac$, we must have $A_k(U)(\xi)\geq0$. Hence, there exists $r>0$ such that  $A_k(U)>0$ on $B(\xi,r)$.\\
Thus, $\Delta_k(U^W\omega_k)\geq0$ in $\mathcal{D}'(B^W(\xi,r))$ i.e. $U^W$ is weakly $\Delta_k$-subharmonic  on $B^W(\xi,r)$ as desired.\\
Now, again, if we follow  the proof of the strong maximum principle in \cite{GalRej3} for  the $W$-invariant $\Delta_k$-subharmonic function $U^W$, then we conclude that
$$
 U=U(\xi)=U(x_0)\quad \text{on}\quad B(\xi,r).
$$
That is $B(\xi,r)\subset \Omega_0$.

\medskip

\noindent
\textbf{Third case:} $E_1\neq \Rac$ and $E_2\neq \emptyset$.  Let $W'\subsetneq W$ be the isotropy group of $\xi$.
Here, we choose $R>0$ under the further following assumption
$$
3R\leq \min_{\widetilde{g}\in \sfrac{W}{W'},\ g\neq id}\|\xi-g\xi\|, \quad \text{with}\quad \sfrac{W}{W'}:=\{\widetilde{g}=gW',\ g\in W\}.
$$
Let $S$ be the $W'$-invariant continuous function  defined on $A$ by $S=\frac{1}{|W'|}\sum_{g'\in W'}g'.U$.\\
$\bullet$ Clearly,  $S$  has a maximum at the point $\xi$ with $S(\xi)=U(\xi)=U(x_0)$. Furthermore, using the $W$-invariance property of $\Delta_k$ as well as  the hypothesis the lemma , we get
\begin{equation}\label{Delta-k S}
\Delta_k(S\omega_k)\geq \frac{1}{|W'|}\sum_{g'\in W'}g'.(\phi(.,v)-\phi(.,u))\omega_k \quad \text{in}\quad \mathcal{D}'(A).
\end{equation}
 $\bullet$  Now, consider the $W$-invariant continuous function $S^W$ defined on
 $$
 B^W(\xi,R):=\cup_{g\in W}B(g\xi,R)=\cup_{\widetilde{g}\in \sfrac{W}{W'}}B(g\xi,R)
 $$
 by
$$
S^W(x):=g.S(x)=S(g^{-1}.x)\quad \text{whenever}\quad x\in B(g\xi,R)\quad \text{and}\  \widetilde{g}\in \sfrac{W}{W'}.
$$
Note that thanks to the previous condition on $R$, the function $S^W$ is well defined. \\
$\bullet$ Let $\widetilde{g_1}=\widetilde{id}$ and $\widetilde{g_2},\dots,\widetilde{g_m}\in \sfrac{W}{W'}$ such that $\textbf{1}_{B^W(\xi,R)}=\sum_{i=1}^m\textbf{1}_{B(g_i\xi,R)}$ and then we can write
$$
 S^W=\sum_{i=1}^mS^W\textbf{1}_{B(g_i\xi,R)}=\sum_{i=1}^mg_i.[S\textbf{1}_{B(\xi,R)}].
$$
$\bullet$ Let $\varphi\in \mathcal{D}(B^W(\xi,R))$ be nonnegative. Following the same idea as in the second case (where we  replace $U$ by $S$ and $U^W$ by $S^W$) and using (\ref{Delta-k S}) we see that we can obtain
\begin{equation}\label{S is subharmonic}
\prs{\Delta_k(S^W\omega_k)}{\varphi}\geq \sum_{i=1}^m\prs{S\omega_k}{A_k(\psi_i)},\quad \text{with}\quad \psi_i=[g_i^{-1}.\varphi]\textbf{1}_{B(\xi,R)}.
\end{equation}
On the other hand, the $W'$-invariance property  of the function $S$ implies that
$$
\forall\ i=1,\dots,m,\ \forall\ \alpha\in E_2,\quad \textstyle\int_{B^W(\xi,R)}S(x)\frac{\psi_i(x)-\psi_i(\sigma_{\alpha}x)}{\prs{\alpha}{x}^2}\omega_k(x)dx=0.
$$
Hence, for every $i=1,\dots,m$ we have
\begin{align*}
\prs{S\omega_k}{A_k(\psi_i)}&=\textstyle\sum_{\alpha\in\Rac_+\setminus E_2}k(\alpha)\|\alpha\|^2\int_{B^W(\xi,R)}S(x)\frac{\psi_i(x)
-\psi_i(\sigma_{\alpha}x)}{\prs{\alpha}{x}^2}\omega_k(x)dx\\
&=\textstyle\sum_{\alpha\in\Rac_+\setminus E_2}k(\alpha)\|\alpha\|^2\int_{B(\xi,R)}\frac{S(x)
-S(\sigma_{\alpha}x)}{\prs{\alpha}{x}^2}[g_i^{-1}.\varphi](x)\omega_k(x)dx\\
&=\prs{A_k(S)\omega_k}{[g_i^{-1}.\varphi]\textbf{1}_{B(\xi,R)}}
\end{align*}
As $E_1\neq\Rac$ and
$$
S(\xi)-S(\sigma_{\alpha}\xi)=\frac{1}{|W'|}\sum_{g'\in W'}\big(U(\xi)-U(\sigma_{g'^{-1}.\alpha}\xi)\big)\geq0,
$$
we deduce that $A_k(S)(\xi)>0$. Consequently, there exists $r>0$ such that $A_k(S)\geq0$ on $B(\xi,r)$. This fact, (\ref{S is subharmonic}), the continuity of $S$ and the Weyl lemma show that $S^W$ is $\Delta_k$-subharmonic on $B^W(\xi,r)$.\\
Now, by the strong maximum principle, we obtain
$$
S=S(\xi)=U(x_0) \quad \text{on}\quad B(\xi,r).
$$
Thus, we get
$$
U=U(\xi)=U(x_0)\quad \text{on}\quad B(\xi,r).
$$
This completes the proof of the lemma.
\fd
\medskip

 The main tool to establish the existence of a solution of the boundary problem (\ref{semilinear-Dirichlet_problem}) is the Schauder  fixed point theorem. In order to apply this theorem,  we will prove  the following intermediate result:
\begin{Prop}\label{Green potential continuity}Let $f$ be a bounded function on $A$ and $G_{k,A}[f]$ be $\Delta_k$-Green potential of $f$ on A given by
\begin{equation}\label{Green potential annulus}
G_{k,A}[f](x):=\int_AG_{k,A}(x,y)f(y)\omega_k(y)dy,\quad x\in A.
\end{equation}
Then $G_{k,A}[f]$ belongs to $\mathcal{C}_0(A)$. Moreover, we have
\begin{equation}\label{Green-Poisson equation}
-\Delta_k\big(G_{k,A}[f]\omega_k\big)=f\omega_k\quad \text{in}\quad \mathcal{D}'(A).
\end{equation}
\end{Prop}
Before proving this result,  we need to show the following  lemma:
\begin{Lem} We have
\begin{equation}\label{Kato condtion}
\lim_{r\rightarrow0}\sup_{x\in A}\eta_{x,r}=0,\quad \text{with}\quad \eta_{x,r}:=\int_{B^W(x,r)}N_k(x,y)\omega_k(y)dy
\end{equation}
and $B^W(x,r):=\cup_{g\in W}B(gx,r)$.
\end{Lem}
\emph{Proof:}
Let $x\in A=A_{\rho,1}$  and $r\in (0,\rho)$.  Since $N_k(x,.)$ is
$\Delta_k$-harmonic on $\R^d\setminus W.x$ and  $\Delta_k$-superharmonic in $\R^d$, by the (super-) mean volume property (\ref{volume mean def})  we deduce that
\begin{align*}
0\leq \eta_{x,r}&\leq \int_{B(0,r+\|x\|)\setminus B(0,\|x\|-r)}N_k(x,y)\omega_k(y)dy\\
&=\int_{B(0,r+\|x\|)}N_k(x,y)\omega_k(y)dy-\int_{B(0,\|x\|-r)}N_k(x,y)\omega_k(y)dy\\
&\leq m_k[B(0,r+\|x\|)]N_k(x,0)-m_k[B(0,\|x\|-r)]N_k(x,0)\\
&= C\ N_k(x,0)\big[(\|x\|+r)^{d+2\gamma}-(\|x\|-r)^{d+2\gamma}\big]\\
&\leq C \frac{\rho^{-2\lambda_k}}{2d_k\lambda_k}\big[(\|x\|+r)^{d+2\gamma}-(\|x\|-r)^{d+2\gamma}\big].
\end{align*}
This shows that $\lim_{r\rightarrow0}\sup_{x\in A}\eta_{x,r}=0$ as desired.
\fd

\medskip

\noindent\emph{Proof of  Proposition} \ref{Green potential continuity}: We can suppose that $f$ is nonnegative.
  Let $\varepsilon>0$. From (\ref{Kato condtion}), there exists $r>0$ such that
\begin{equation}\label{Kato condtion continuity proof}
\forall\ x\in A,\quad \eta_{x,2r}<\varepsilon.
\end{equation}
$\bullet$ First, we will prove that $G_{k,A}[f](x)\longrightarrow0$ when $x$ tends to $\partial{A}$.
Let $x\in A$. By (\ref{Kato condtion continuity proof}) we have
\begin{align*}
G_{k,A}[f](x)&=\int_{A}G_{k,A}(x,y)f(y)\omega_k(y)dy\\
&= \int_{A\cap B^W(x,r)}G_{k,A}(x,y)f(y)\omega_k(y)dy+\int_{A\setminus{B}^W(x,r)}G_{k,A}(x,y)f(y)\omega_k(y)dy\\
&\leq \|f\|_{\infty}\eta_{x,r}+\|f\|_{\infty}\int_{A\setminus{B}^W(x,r)}G_{k,A}(x,y)\omega_k(y)dy\\
&\leq \varepsilon\|f\|_{\infty}+\|f\|_{\infty}\int_{A\setminus{B}^W(x,r)}G_{k,A}(x,y)\omega_k(y)dy.
\end{align*}
Since for every $z\in \text{supp}\ \mu_y\subset Co(y)$, we can write
$$
\|x\|^2+\|y\|^2-2\prs{x}{z}=\textstyle\sum_{g\in W}\lambda_g(z)\|x-gy\|^2,
$$
 with $\lambda_g(z)\geq0$ and $\sum_{g\in W}\lambda_g(z)=1$, we deduce that
$$
\forall\ y\in A\setminus{B}^W(x,r),\quad 0\leq G_{k,A}(x,y)\leq N_k(x,y)\leq\frac{r^{-2\lambda_k}}{2d_k\lambda_k}.
$$
Hence, we can apply the Lebesgue dominated convergence theorem to obtain
$$
\lim_{x\rightarrow\xi\in\partial{A}}\int_{A\setminus{B}^W(x,r)}G_{k,A}(x,y)\omega_k(y)dy=0.
$$
\medskip

\noindent $\bullet$ Now, we will prove that $G_{k,A}[f]$ is continuous on $A$.  Fix $x_0\in A$ and  assume that $\overline{B}(x_0,2r)\subset A$.
Since $f$ is bounded, it is enough to prove that
\begin{equation}\label{Ascoli}
\lim_{x\rightarrow x_0}\int_{A}|G_{k,A}(x,y)-G_{k,A}(x_0,y)|\omega_k(y)dy=0
\end{equation}
For any $x\in B(x_0,r)$, we have
\begin{align*}
\int_{A}|G_{k,A}(x,y)-G_{k,A}(x_0,y)|\omega_k(y)dy
&\leq \int_{B^W(x_0,r)}|G_{k,A}(x,y)-G_{k,A}(x_0,y)|\omega_k(y)dy\\
&+\int_{A\setminus B^W(x_0,r)}|G_{k,A}(x,y)-G_{k,A}(x_0,y)|\omega_k(y)dy\\
&=I_1(x,x_0)+I_2(x,x_0).
\end{align*}
As $B^W(x_0,r)\subset B^W(x,2r)$, by (\ref{Kato condtion continuity proof}) we have
\begin{align*}
I_1(x,x_0)&\leq \int_{B^W(x_0,r)}N_k(x,y)\omega_k(y)dy+\int_{B^W(x_0,r)}N_k(x_0,y)\omega_k(y)dy\\
&\leq\eta_{x,2r}+\eta_{x_0,r}\leq 2\varepsilon.
\end{align*}
In addition, by item 8) in Proposition \ref{Green function properties},  we know that the function $(x,y)\longmapsto G_{k,A}(x,y)$ is continuous on the compact set $\overline{B}^W(x_0,r)\times \big(\overline{A}\setminus B^W(x_0,r)\big)$. Thus, there exists $\theta>0$ such that for every $x\in B(x_0,\theta)$ and every $y\in A\setminus B^W(x_0,r)$, we have
$$
 |G_{k,A}(x,y)-G_{k,A}(x_0,y)|\leq \varepsilon.
$$
This implies that
\begin{align*}
\forall\ x\in B(x_0,\theta),\quad  I_2(x,x_0)\leq \varepsilon \int_{A}\omega_k(y)dy.
\end{align*}
Finally, we conclude that $G_{k,A}[f]\in \mathcal{C}_0(A)$.\\
\noindent$\bullet$ Let $\varphi\in \mathcal{D}(A)$. Using Fubini's theorem, the symmetry property of the Green function and (\ref{Green-Poisson equation1}) we get
\begin{align*}
-\prs{\Delta_k\big(G_{k,A}[f]\omega_k\big)}{\varphi}&=-\int_{A}f(y)\prs{\Delta_k\big(G_{k,A}(.,y)\omega_k\big)}{\varphi}\omega_k(y)dy \\
&=\int_{A}f(y)\varphi(y)\omega_k(y)dy.
\end{align*}
This completes the proof.
\fd

\bigskip

\noindent \emph{Proof of  Theorem} \ref{thm semilinear problem}:  Fix $f\in \mathcal{C}^+(\partial{A})$ and
\begin{align*}
 c_1&:=\inf_{x\in{\overline{A}}}\big(P_{k,A}[f](x)-G_{k,A}[\phi(.,c_2)](x)\big)\\
 &=\inf_{x\in{\overline{A}}}\big(P_{k,A}[f](x)-\phi_2(c_2)G_{k,A}[\phi_1](x)\big),\quad\text{with} \quad c_2:=\max_{\overline{A}}P_{k,A}[f].
\end{align*}
Let us consider the bounded, closed and convex set
$$
\mathcal{M}:=\{u\in \mathcal{C}(\overline{A}):\quad c_1\leq u\leq c_2\}
$$
endowed the uniform topology and the map $T:\mathcal{C}(\overline{A})\longrightarrow \mathcal{C}(\overline{A})$ defined by
$$
T(u):=P_{k,A}[f]-G_{k,A}\big(\phi(.,u)\big).
$$
Note that since $\phi(x,u(x))=\phi_1(x)\phi_2(u(x))$ is bounded, by Proposition \ref{Green potential continuity}, $G_{k,A}\big(\phi(.,u)\big)\in \mathcal{C}_0(A)$ and
then $T$ is well defined. Moreover, as $\phi_2$ is nondecreasing, for every $u\in \mathcal{M}$ and every $x\in \overline{A}$, we have
\begin{equation}\label{pointwise bounded}
c_1\leq P_{k,A}[f](x)-G_{k,A}[\phi(.,c_2)](x)\leq T(u)(x)\leq P_{k,A}[f](x)\leq c_2.
\end{equation}
Hence, we have $T(\mathcal{M})\subset \mathcal{M}$.\\
Now, we want to establish that $T$ has a unique fixed point in $\mathcal{M}$ by using the Schauder theorem. \\
$\bullet$ Firstly, we will prove that $T(\mathcal{M})$ is relatively compact.  For this, we will use the Arzel\`{a}-Ascoli theorem.
From (\ref{pointwise bounded}), $T(\mathcal{M})$ is pointwise bounded. \\
Let $x_0\in A$. For every $u\in \mathcal{M}$  we have
\begin{align*}
|T(u)(x)-T(u)(x_0)|&\leq \big|P_{k,A}[f](x)-P_{k,A}[f](x_0)\big|+\big|G_{k,A}[\phi(.,u)](x)-G_{k,A}[\phi(.,u)](x_0)\big|\\
&\leq \big|P_{k,A}[f](x)-P_{k,A}[f](x_0)\big|\\
&+\int_A\big|G_{k,A}(x,y)-G_{k,A}(x_0,y)\big|\phi_1(y)\phi_2(u(y))\omega_k(y)dy\\
&\leq \big|P_{k,A}[f](x)-P_{k,A}[f](x_0)\big|\\
&+\phi_2(c_2)\|\phi_1\|_{\infty}\int_A\big|G_{k,A}(x,y)-G_{k,A}(x_0,y)\big|\omega_k(y)dy.
\end{align*}
Therefore, from (\ref{Ascoli}) and the continuity of the function $P_{k,A}[f]$, we conclude that $T(\mathcal{M})$ is equicontinuous. Finally,
$T(\mathcal{M})$ is relatively compact as desired.  \\
$\bullet$ Secondly, we will prove that $T:\mathcal{M}\longrightarrow \mathcal{M}$ is continuous. Let then $(u_n)$ be sequence in $\mathcal{M}$ which
converges uniformly to $u\in \mathcal{M}$. We have
 \begin{align*}
|T(u_n)(x)-T(u)(x)|&\leq\int_AG_{k,A}(x,y)\phi_1(y)\big|\phi_2(u_n(y))-\phi_2(u(y))\big|\omega_k(y)dy.
\end{align*}
But
$$
0\leq G_{k,A}(x,y)\phi_1(y)\big|\phi_2(u_n(y))-\phi_2(u(y))\big|\leq 2\phi_2(c_2)\|\phi_1\|_{\infty}G_{k,A}(x,y).
$$
Thus, we can use the Lebesgue dominated convergence theorem to obtain that $T(u_n)\longrightarrow T(u)$ pointwise. Hence, by equicontinuity, we get the uniform convergence.\\
Consequently, there exists $u\in \mathcal{M}$ such that
$$
u=T(u)=P_{k,A}[f]-G_{k,A}\big(\phi(.,u)\big).
$$
Finally, note that from the properties of $P_{k,A}$ as well as (\ref{Green-Poisson equation}), $u$ is a solution of (\ref{semilinear-Dirichlet_problem}).
This finishes  the proof of the theorem.

\fd

\begin{center}
   \emph{Acknowledgement}
\end{center}
It is a pleasure to thank the referee for the valuable suggestions which
improved the presentation of the paper.



\begin{thebibliography}{XX}
\scriptsize
\bibitem{Axler} S. Axler, P. Bourdon and W. Ramey. \emph{Harmonic Function Theory}. Springer Verlag , Second edition (2001).
\bibitem{Ben Chrouda}M. Ben Chrouda.  \emph{On the Dirichlet problem associated with the Dunkl Laplacian}. Ann. Polon. Math.117(1), (2016), 79-87.
\bibitem{Ben Chrouda 2}M. Ben Chrouda, K. El Mabrouk and K. Hassine. Boundary value problem for the Dunkl Laplacian. Accepted in the Journal of Prob. and Math. Stat.
\bibitem{Diejen-Vinet} J. F. van Diejen and L. Vinet. \emph{Calogero-Sutherland-Moser Models}. Springer-Verlag, CRM Series in Mathematical Physics (2000).
\bibitem{DaiXu} F. Dai and Y. Xu. \emph{Approximation Theory and Harmonic Analysis on Spheres and Balls}. Springer, (2013).
\bibitem{Dunkl1} C. F. Dunkl. \emph{Differential-difference operators associated to reflection groups}.
Trans. Amer. Math. Soc., 311, (1989), 167-183.
\bibitem{Dunkl2}C. F. Dunkl. \emph{Integral kernels with reflection group invariance}. Canad. J. Math., 43, (1991), 123-183.
\bibitem{DunklXu} C. F. Dunkl and Y. Xu. \emph{Orthogonal Polynomials of Several variables}. Cambridge Univ. Press (2001).
\bibitem{El Kamel} J. El Kamel and C. Yacoub. \emph{Poisson integrals and Kelvin transform associated to Dunkl-Laplacian operator}. Global Journal of Pure and Applied Math. (2007), Vol. 3, Issue 5, p. 351.
\bibitem{Etingof CMS} P. Etingof. \emph{Calogero Moser systems and representation theory}. Z\"{u}rich Lectures in Advanced Mathematics,
European Mathematical Society (EMS), Z\"{u}rich, (2007).
\bibitem{GalRej} L. Gallardo and C. Rejeb. \emph{A new mean value property for harmonic functions relative to the Dunkl-Laplacian operator and applications}.
Trans.  Amer. Math. Soc., Vol. 368, Number 5, May 2016, p.3727-3753.
\bibitem{GalRej2} L. Gallardo and C. Rejeb. \emph{Radial mollifiers, mean value operators and harmonic functions in Dunkl theory}.  J.
Math. Anal. Appl. (2017), Volume 447, Issue 2, 1142-1162.
\bibitem{GalRej3} L. Gallardo and C. Rejeb. \emph{Newtonian Potentials and subharmonic functions associated to root systems}. J. Potential Anal, 47 (2017), 369-400.
\bibitem{GalRejSifi} L. Gallardo, C. Rejeb and M. Sifi. \emph{Riesz potentials of Radon measures associated to reflection groups}. Adv. in Pure and Appl. Math, Vol 9, (2018), 109-130.
\bibitem{Graczyk-Luks-Rosler} P. Graczyk, T. Luks and M. R\"{o}sler. \emph{On the Green Function and Poisson Integrals of the Dunkl Laplacian}. J.  Potential Anal, Volume 48, Issue 3 (2018),  337-360.
\bibitem{Grossi} M. Grossi and D. Vujadinovi\'{c}. \emph{On the Green Function of the Annulus}. Anal. Theory Appl., Vol. 32, No. 1 (2016), 52-64.
\bibitem{Kods} K. Hassine.  Mean value property of $\Delta_k$-harmonic functions on $W$-invariant open sets. Afr. Mat. 27(7), (2016), 1275-1286.
\bibitem{humph} J. E. Humphreys. \emph{Reflection groups and Coxeter groups}. Cambridge Studies in Advanced Mathematics 29, Cambridge University Press, (1990).
\bibitem{Maslyouss} M. Maslouhi and E. H. Youssfi. \emph{Harmonic functions associated to Dunkl operators}. Monatsh. Math. 152 (2007), 337-345.
\bibitem{MajTrim}H. Mejjaoli and  K. Trim\`{e}che . \emph{On a mean value property associated with the Dunkl Laplacian operator and
applications}. Integ. Transf. and Spec. Funct., 12(3), (2001), 279-302.
\bibitem{Rosler3}M. R\"{o}sler. \emph{Generalized Hermite polynomials and the heat equation for Dunkl operators}. Comm. Math. Phys, 192, (1998), 519-542.
\bibitem{Rosler-Voit}M. R\"{o}sler and M. Voit. \emph{Markov processes related with Dunkl operators}. Adv. in Appl. Math. 21 (1998), 575-643.
\bibitem{Rosler1}M. R\"{o}sler. \emph{Positivity of Dunkl's intertwining operator}. Duke Math. J., 98, (1999), 445-463.
\bibitem{Rosler4}M. R\"{o}sler. \emph{Dunkl Operators: Theory and Applications}. Lecture Notes in Math., vol.1817, Springer Verlag (2003), 93-136.
\bibitem{Szego} G. Szeg\"{o}. \emph{Orthogonal Polynomials}. Amer. Math. Soc., Providence, RI, Fourth edition, (1975).
\bibitem{Trimeche1}K. Trim\`{e}che. \emph{The Dunkl intertwining operator on spaces of functions and distributions and integral representation
of its dual}. Integ. Transf. and Spec. Funct., 12(4), (2001), 394-374.
\bibitem{Xu97} Y. Xu. \emph{Integration of the intertwining operator for h-harmonic polynomials associated to reflection groups}. Proc. Amer. Math. Soc. 125 (1997), 2963-2973.
\end{thebibliography}
\end{document}